\documentclass{elsart}

\usepackage{amssymb,amsbsy}
\usepackage[mathcal]{euscript}

\def\u{\stackrel{\mbox{\tiny $\rightarrow$}}{\boldsymbol u}}
\def\vv{\stackrel{\mbox{\tiny $\rightarrow$}}{\boldsymbol v}}
\def\dd{\delta}
\def\bphi{\bar{\phi}}
\def\bchi{\bar{\chi}}
\def\sphi{\phi^*}
\def\schi{\chi^*}
\def\rmg{{\rm g}} 
\def\Gl{{\mathcal G}}
\def\nb{\nabla}
\def\bnb{\bar{{\nabla}}}
\def\bR{\bar{R}}

\def\U{{\mathcal U}}

\def\G{\boldsymbol{\rmg}}

\def\bP{\boldsymbol{P}}
\def\bPi{\boldsymbol{\Pi}}
\def\bfeta{\boldsymbol{\eta}}
\def\bz{\bar{\boldsymbol{\theta}}}

\def\z{\theta}

\def\fr{\frac}
\def\d{\partial}

\def\g{\gamma}
\def\Chr{\Gamma}

\def\bbg{\bar{\gamma}}
\def\vecX{{\mathfrak X}}
\def\xiv{\stackrel{\mbox{\tiny $\rightarrow$}}{\boldsymbol\xi}}
\def\lie{\pounds_{\xiv}}

\def\be{\begin{equation}}
\def\ee{\end{equation}}
\def\bnr{\begin{eqnarray*}}
\def\enr{\end{eqnarray*}}
\def\bea{\begin{eqnarray}}
\def\eea{\end{eqnarray}}

\begin{document}

\begin{frontmatter}



\title{Bi-conformal vector fields and 
the local geometric characterization of conformally 
separable pseudo-Riemannian manifolds I}

\author{Alfonso Garc\'{\i}a-Parrado G\'omez-Lobo}

\address{Departamento de F\'{\i}sica Te\'orica, 
Universidad del Pa\'{\i}s Vasco. Apartado 644, 48080 Bilbao 
(Spain)}

\ead{wtbgagoa@lg.ehu.es}

\begin{abstract}
This is the first of two companion papers in which 
a thorough study of the normal form and the first integrability 
conditions arising from {\em bi-conformal vector fields} is presented. 
These new symmetry transformations were introduced in 
{\em Class. Quantum Grav.} \textbf{21}, 2153-2177
and some of their basic properties were addressed there. Bi-conformal vector fields 
are defined on a pseudo-Riemannian manifold $V$
through the differential conditions $\lie P_{ab}=\phi P_{ab}$ and $\lie\Pi_{ab}=\chi\Pi_{ab}$ where 
$P_{ab}$ and $\Pi_{ab}$ are orthogonal and complementary projectors with respect to the 
metric tensor $\rmg_{ab}$. In our calculations a new affine connection ({\em bi-conformal
connection}) arises quite naturally and this connection enables us to find a local 
characterization of {\em conformally separable} pseudo-Riemannian manifolds (also called double 
twisted products) in terms of the vanishing of a rank three tensor $T_{abc}$. 
Similar local characterizations are found
for the most important particular cases such as (double) warped products, 
twisted products and conformally reducible spaces.
\end{abstract}

\end{frontmatter}

\section{Introduction}
\label{intro}
The research of symmetry transformations in Differential Geometry and General Re\-la\-ti\-vity has been an 
important subject during 
the years.
Here by symmetries we mean a group of transformations of a given pseudo-Riemannian manifold complying with certain geometric 
property. By far the most studied symmetries are isometries and conformal transformations which are defined 
through the conditions
\be
\lie\rmg_{ab}=0,\ \ \lie\rmg_{ab}=2\phi\rmg_{ab},
\label{n-zero}
\ee
where $\rmg_{ab}$ is the metric tensor of the manifold, $\xiv$ is the {\em infinitesimal generator} of the transformation and
$\phi$ is a function which we will call {\em gauge} of the symmetry (this terminology was first 
employed in \cite{SERGI} and it will be explained later). 
Infinitesimal generators of these symmetries are known as Killing vectors and conformal Killing vectors respectively. 
As is very simple to check they are a Lie algebra with respect to the Lie bracket of vector fields 
and the transformations generated by these vector fields give rise to subgroups of the diffeomorphism group.  

Important questions are the possible dimensions of these Lie algebras and the geometric characterizations of spaces admitting the
symmetry. The general answer to these questions can in principle be obtained by solving the differential conditions 
written above although for general enough cases the explicit evaluation of such solutions gets too complex and other methods are 
required. Notwithstanding these difficulties, we can obtain easily from the differential conditions the cases in which the Lie 
algebras are finite dimensional, the greatest dimension of these Lie algebras and geometric characterizations of the 
spaces admitting these Lie algebras as solutions. This is done by finding the {\em normal form} of the above equations 
(if such form exists) and the complete integrability conditions coming from this set of equations.
In this way we deduce that isometries are always finite dimensional whereas 
conformal motions are finite dimensional iff the space dimension is greater or equal than three. 
The spaces in which the greatest 
dimension is achieved are constant curvature and conformally flat spaces respectively
and as is very well known they are characterized by the geometric conditions
\bnr
R^a_{\ bcd}&&=\fr{R}{n(n-1)}(\dd^a_{\ c}\rmg_{bd}-\dd^a_{\ d}\rmg_{bc}),\ 
\mbox{(constant curvature)},\\ 
C^a_{\ bcd}&&=0,\ n>3,\ \mbox{(conformally flat)}
\label{weyl}
\enr
where $n$ is the dimension of the manifold, $R^a_{\ bcd}$ is the curvature tensor,
 $R$ the scalar curvature and $C^a_{\ bcd}$ 
the Weyl tensor\footnote{In the case of dimension three manifolds the Weyl tensor 
is replaced by a three rank tensor called the Cotton-York (or Schouten) tensor.}

The procedure followed for isometries and conformal motions is carried over 
to other symmetries such as 
linear and affine collineations and conformal collineations 
(see \cite{YANO,LIBROHALL} for a very good 
account of this). 
However, little research has been done for 
symmetries different from these mostly because the cases under 
consideration were infinite dimensional {\em generically}. 
This means that it is not possible to obtain a normal set of 
equations out of the differential conditions (see section \ref{normal-dimension}) 
which greatly complicates matters. Mathematicians have developed an alternative
view toward this issue in the theory of $G$-{\em structures} (see e.g. \cite{REINHART} 
for a thorough description of this).

In reference \cite{BI-CONFORMAL} we put forward a new symmetry 
transformation for general pseudo-Riemannian manifolds. 
Infinitesimal generators of these symmetries (bi-con\-for\-mal vector fields) 
fulfill the differential conditions 
\be
\lie P_{ab}=\phi P_{ab},\ \ \lie\Pi_{ab}=\chi\Pi_{ab},
\label{n-uno}
\ee
where $P_{ab}$ and $\Pi_{ab}$ are orthogonal and complementary
 projectors with respect to the metric tensor $\rmg_{ab}$ and 
$\phi$, $\chi$ are the gauges of the symmetry. These are 
functions which, as happened in the conformal case, depend on the 
vector field $\xiv$ so a solution of (\ref{n-uno}) is formed by 
$\xiv$ itself and the gauges $\phi$ and $\chi$ (we will 
usually omit the dependence on $\xiv$ in the gauges).
The finite transformations generated by bi-conformal vector
 fields are called bi-conformal transformations. 
In a sense, these symmetries can be regarded as conformal 
transformations with respect to both $P_{ab}$ and $\Pi_{ab}$ 
so we can expect that some properties of bi-conformal vector 
fields will resemble 
those of conformal transformations. In \cite{BI-CONFORMAL} it 
was shown that bi-conformal vector fields comprise a Lie algebra 
under the Lie bracket and that this algebra is finite dimensional
if none of the projectors has algebraic rank one or two being the 
greatest dimension
$$
N=\fr{1}{2}p(p+1)+\fr{1}{2}(n-p)(n-p+1),
$$
with $p$ the algebraic rank of one of the projectors. We provided also explicit examples in which this dimension is achieved,
namely bi-conformally flat spaces which in local coordinates $x=\{x^1,\dots,x^n\}$ look like 
$(\alpha,\beta=1,\dots p,\ A,B=p+1,\dots n)$ 
\be
ds^2=\Xi_1(x)\eta_{\alpha\beta}dx^{\alpha}dx^{\beta}+\Xi_2(x)\eta_{AB}dx^Adx^B,\ \Xi_1,\ \Xi_2\in C^3,
\label{conformal-1}
\ee
where $\eta_{\alpha\beta}$, $\eta_{AB}$ are flat metrics depending only on 
the coordinates $x^{\alpha}$ and $x^A$ respectively.
That
these spaces play the same role for bi-conformal vector fields as conformally flat spaces or spaces of 
constant curvature do for the classical symmetries will be a result 
of the analysis started in this paper.
One can also find a geometric characterization for 
bi-conformally flat spaces similar to those 
of the spaces of constant curvature or conformally flat spaces 
stated in (\ref{weyl}) 
(full details of this are contained in \cite{GRQC}). In the scheme developed 
in \cite{BI-CONFORMAL} this
sort of characterization could not be extracted due to the complexity of 
the calculations and it had to be postponed. 

In this paper we perform the full calculation of the normal form 
for equations (\ref{n-uno}). This normal form is already present in our
previous work but it turned out to be rather messy and relevant geometric
information could not be obtained. 
This was so because all these calculations were done using the covariant 
derivatives arising from the metric
connection which is not adapted to the calculations. Here we show that the 
definition of a new symmetric connection ({\em bi-conformal connection}) 
greatly simplifies the calculations making it possible to get a simpler form 
for the normal system. Due to the
great amount of algebra required to work out the complete integrability 
conditions associated to the normal form we 
have placed its analysis together with the geometric characterization 
of the maximal spaces in a subsequent 
paper (a complete version of all our results can be found in \cite{GRQC}).  

The bi-conformal connection bears an interesting geometric interpretation 
if we work on {\em conformally separable} pseudo-Riemannian manifolds. 
These are defined as those manifolds
which in a local coordinate system the metric tensor takes the form 
(the conventions are the same as in (\ref{conformal-1}))
$$
ds^2=\Xi_1(x)g_{\alpha\beta}dx^{\alpha}dx^{\beta}+\Xi_2(x)G_{AB}dx^Adx^B,
$$
where $g_{\alpha\beta}$ and $G_{AB}$ only depend on the coordinates labelled by their index components.
The bi-conformal connection is naturally adapted to this decomposition and its use permits us to give a new simple 
geometric characterization of these spaces in terms of the 
vanishing of a certain rank-three tensor $T_{abc}$. The most known 
cases of conformally separable pseudo-Riemannian manifolds are warped products,
double warped products, twisted products and conformally reducible 
spaces (see definition 
\ref{classification} for a precise account of each case)
and we can easily derive with our techniques local 
geometric characterizations for these spaces (theorem \ref{theorem-classification}).
   
The outline of the paper is as follows: section \ref{definitions} introduces the basic notation and definitions. 
 In section \ref{connection} we define a new symmetric connection (bi-conformal connection) and we set its main properties.
Section \ref{normal-dimension} presents the calculation of the normal form associated to (\ref{n-uno}) and the calculation
of the maximum dimension of any finite dimensional Lie algebra of bi-conformal vector fields is carried out. 
In section \ref{sufficient} we use the bi-conformal connection to supply a local geometric 
characterization of conformally separable pseudo-Riemannian manifolds
 and their principal subcases. Finally in section \ref{examples} 
we show in explicit examples how to use this geometric characterization 
and we hint how these conditions may be extended 
to more general pseudo-Riemannian manifolds. 
Appendix A collects basic identities relating
the Lie derivative and the covariant derivative.  

\section{Bi-conformal vector fields and bi-conformal transformations}
\label{definitions}
Let us start by setting our notation and conventions for the paper.
We work on a differentiable manifold $V$ in which a $C^{\infty}$
metric $\rmg_{ab}$ of arbitrary signature has been defined (pseudo-Riemannian manifold). 
Vectors and vector fields are denoted with arrowed characters 
$\u$, $\vv$
 (we leave to the context the distinction between each of these 
entities) when expressed in coordinate-free notation whereas 
1-forms are written in bold
characters ${\boldsymbol u}$. Sometimes this same notation will be employed for other higher rank objects
such as contravariant and covariant tensors. 
 Indexes of tensors are represented by lowercase Latin characters $a$, $b,\dots $ 
and the metric $\rmg_{ab}$ or its inverse $\rmg^{ab}$ are used to respectively raise 
or lower indexes. Rounded and square brackets are used for symmetrization and antisymmetrization respectively
and whenever a group of indexes is enclosed between strokes they are excluded from the symmetrization 
or antisymmetrization operation.
 Partial derivatives with respect to local coordinates are $\d_a\equiv\d/\d x^a$.
The Levi-Civita connection associated to $\rmg_{ab}$ is $\g^{a}_{\ bc}$ (Ricci rotation coefficients) reserving
the symbol $\Gamma^a_{\ bc}$ only for the Christoffel symbols, namely, the connection components
in a natural basis. The covariant derivative and the Riemann tensor constructed from 
this connection are denoted by $\nb$ and $R^a_{\ bcd}$ respectively being our convention for the 
Riemann tensor
\be
R^a_{\ bcd}\equiv\d_c\Chr^a_{\ db}-\d_d\Chr^a_{\ cb}+\Chr^a_{\ rc}\Chr^r_{\ db}-\Chr^a_{\ rd}\Chr^r_{\ cb}.
\label{convention}
\ee
Under this convention the Ricci identity becomes 
$$
\nb_b\nb_cu^a-\nb_c\nb_bu^a=R^a_{\ rbc}u^r,\ \ \nb_b\nb_cu_a-\nb_c\nb_bu_a=-R^r_{\ abc}u_r.
$$ 
All the above relations are still valid for a non-metric symmetric connection.

The set of smooth vector fields of the manifold
$V$ is denoted by $\vecX(V)$. This is an infinite dimensional Lie algebra  which is sometimes regarded as the Lie
algebra of the group of diffeomorphisms of the manifold $V$. Finally the Lie derivative operator 
with respect to a vector field $\xiv$ is $\lie$.

One of the main subjects of this paper is the study of bi-conformal vector fields whose definition 
given in \cite{BI-CONFORMAL} we reproduce here.
\begin{defn}
A smooth vector field $\xiv$ on $V$ is said to be a 
{\em bi-conformal vector field} if it fulfills the condition
\be 
\lie P_{ab}=\phi P_{ab},\ \ \ 
\lie\Pi_{ab}=\chi\Pi_{ab},
\label{bi-conformal}
\ee
for some functions $\phi$, $\chi\in C^{\infty}(V)$. 
\label{Bi-conformal}
\end{defn}
$P_{ab}$ and $\Pi_{ab}$ are smooth sections of the tensor
bundle $T^0_2(V)$ such that at each point $x\in V$ they form a pair of orthogonal and 
complementary projectors with respect to the metric tensor $\rmg_{ab}|_x$. This leads to 
$$
P_{ab}=P_{ba},\ \Pi_{ab}=\Pi_{ba},\ 
P_{ab}+\Pi_{ab}=\rmg_{ab},\ P_{ap}P^{p}_{\ b}=P_{ab},\ \Pi_{ap}\Pi^{p}_{\ b}=\Pi_{ab},\  
P_{ap}\Pi^p_{\ b}=0.
$$ 
Equation (\ref{bi-conformal}) can be re-written in a 
number of equivalent ways as shown next.
To that end we define the tensor $S_{ab}$ in terms of the 
projectors $P_{ab}$ and $\Pi_{ab}$ by 
\be
\hspace{-1cm}S_{ab}\equiv P_{ab}-\Pi_{ab}, \Rightarrow P_{ab}=\fr{1}{2}(\rmg_{ab}+S_{ab}),\ \ 
\Pi_{ab}=\fr{1}{2}(\rmg_{ab}-S_{ab})\Rightarrow S_{ap}S^p_{\ b}=\rmg_{ab}. 
\label{square-root}
\ee
The last property of this set means that $S_{ab}$ is a square root of 
the metric tensor.  It is not difficult to prove that 
the endomorphism $S^a_{\ b}$ can be always diagonalized and the only possible 
eigenvalues are $+1$ and $-1$ being the associated eigenspaces the subspaces 
upon which the projectors $P^a_{\ b}$ and $\Pi^a_{\ b}$ project. 
Other interesting point is that a square root is always the {\em superenergy} of 
a simple form (see \cite{BI-CONFORMAL} for further details).

In terms of the square root $S_{ab}$
the conditions (\ref{bi-conformal}) take the equivalent form
\be
\hspace{-1cm}\lie\rmg_{ab}=\alpha\rmg_{ab}+\beta S_{ab},\ \ 
\lie S_{ab}=\alpha S_{ab}+\beta\rmg_{ab},\ \ 
\alpha=\fr{1}{2}(\phi+\chi),\ \ \beta=\fr{1}{2}(\phi-\chi).
\label{second}
\ee
From equation (\ref{square-root}) 
we deduce that both projectors are fixed by the square root $S_{ab}$
so we can use the latter instead of the projectors when working with a given set of 
bi-conformal vector fields. 
Following \cite{BI-CONFORMAL} the set of bi-conformal vector fields possessing $S_{ab}$ as the 
associated square root will be denoted by $\Gl(S)$. 
In this paper only expressions involving $P_{ab}$ and $\Pi_{ab}$ will be used in our calculations.
A very important property of $\Gl(S)$ 
is that it forms a Lie subalgebra of $\vecX(V)$
 (proposition 5.2 of \cite{BI-CONFORMAL}) which can be finite or infinite dimensional.
Conditions upon the tensor $S_{ab}$ (or equivalently the projectors) for this Lie algebra to be finite dimensional 
were given in \cite{BI-CONFORMAL} and they will be re-derived in section \ref{normal-dimension} in a more efficient 
way. Observe that the functions $\phi$ and 
$\chi$ appearing in definition \ref{Bi-conformal} (or $\alpha$ and $\beta$)
do depend on the bi-conformal vector field $\xiv$ (this dependence can be
dropped if we work with a single bi-conformal vector field but 
it should be added when working with Lie algebras of 
bi-conformal vector fields). In the latter case $\phi$ 
and $\chi$ ($\alpha$ and $\beta$) are called {\em gauge functions}
(see \cite{SERGI} for an explanation of this terminology).   

The next set of relations comes straight away from (\ref{bi-conformal})
\be
\lie P^a_{\ b}=\lie \Pi^a_{\ b}=0,\ \ \lie P^{ab}=-\phi P^{ab},\ \lie\Pi^{ab}=-\chi\Pi^{ab}.
\ee
Here the last pair of equations are equivalent to (\ref{bi-conformal}).

\section{The bi-conformal connection}
\label{connection}
As we have commented in the introduction
the Levi-Civita connection $\gamma^a_{\ bc}$ 
is not suitable to study the normal form and the integrability 
conditions coming from the differential condition (\ref{bi-conformal}) as they result in 
rather cumbersome expressions. In order to proceed further in our study we are going to 
show next that the definition of a new symmetric connection greatly simplifies 
the normal form calculated in \cite{BI-CONFORMAL} and what is more, it will enable us 
to work out thoroughly the 
complete integrability conditions arising from this normal form in a subsequent work.

To start with we recall some identities satisfied by any bi-conformal vector field $\xiv$
which were obtained in \cite{BI-CONFORMAL}. These identities are in fact linear combinations of 
the first covariant
derivative of (\ref{bi-conformal}) and we also indicate briefly 
how they are obtained as this information
will be needed later. Using equation (\ref{lie-connection}) we easily obtain the Lie derivative 
of the metric connection $\g^a_{\ bc}$ ($\phi_b\equiv\d_b\phi$, $\chi_b\equiv\d_b\chi$)
\be
\hspace{-1cm}
\lie\g^a_{bc}=\fr{1}{2}(\phi_bP^a_{\ c}+\phi_cP^a_{\ b}-\phi^aP_{bc}+
\chi_b\Pi^a_{\ c}+
\chi_c\Pi^a_{\ b}-\chi^a\Pi_{cb}+(\phi-\chi)M^a_{\ bc}),
\label{3}
\ee
where the tensor $M_{abc}$ is defined by
\be
M_{abc}\equiv\nb_bP_{ac}+\nb_cP_{ab}-\nb_aP_{bc}.
\ee
 
The Lie derivative of $M_{abc}$ can be worked out by means of (\ref{lie-conmmutation}) 
getting
\bea
\lie M_{abc}=\phi M_{abc}+(\chi-\phi)P_{ap}M^{p}_{\ bc}-       
P_{bc}\Pi_{ap}\phi^p+\Pi_{cb}P_{ap}\chi^p=\nonumber\\
=\chi M_{abc}+(\phi-\chi)\Pi_{ap}M^p_{\ bc}-       
P_{bc}\Pi_{ap}\phi^p+\Pi_{cb}P_{ap}\chi^p\label{9},
\eea
from which, projecting down with $P^{bc}$ and $\Pi^{bc}$, we deduce 
\be
\lie E_a=-p\Pi_{ap}\phi^p,\ \ \ \lie W_a=(p-n)P_{ap}\chi^p,
\label{18}
\ee
with the definitions 
$$
E_a\equiv M_{acb}P^{cb},\ \ W_a\equiv-M_{acb}\Pi^{cb},\ p=P^a_{\ a}.
$$
The following algebraic properties of the tensors $E_a$ and $W_a$ are useful
\be
\Pi_{ac}E^c=E_a,\  P_{ac}W^c=W_a,\ 
0=P^{ab}E_b=\Pi^{ab}W_b.
\label{vwe}
\ee
Now, we substitute (\ref{18}) into (\ref{9}) yielding
\bea
\lie\left(M_{acb}-\fr{1}{p}E_aP_{bc}+\fr{1}{n-p}W_a\Pi_{cb}\right)=
\phi\left(\Pi_{ap}M^p_{\ cb}-\fr{E_aP_{bc}}{p}\right)+\nonumber\\
\hspace{15mm}+\chi\left(P_{ap}M^p_{\ bc}+\fr{\Pi_{cb}W_a}{n-p}\right).
\label{integral-constraint}
\eea 
This equation can be written in a more compact form
\be
\lie T_{abc}=(\phi\Pi_{ar}+\chi P_{ar})T^r_{\ bc}=\phi B_{abc}+\chi A_{abc},
\label{compact}
\ee
where the definitions of the tensors $T_{abc}$, $A_{abc}$, $B_{abc}$ are
\bea
T_{abc}\equiv M_{abc}+\fr{1}{n-p}W_a\Pi_{bc}-\fr{1}{p}E_aP_{bc},\label{tensor-t}\\ 
A_{abc}\equiv P_{a}^{\ d}T_{dbc}=P_a^{\ d}M_{dcb}+\fr{1}{n-p}W_a\Pi_{cb},\label{splitting-1}\\
B_{abc}\equiv\Pi_{a}^{\ d}T_{dbc}=\Pi_{a}^{\ d}M_{dcb}-\fr{1}{p}E_aP_{cb}.\label{splitting-2}
\eea
Using equation (\ref{compact}) we can calculate the Lie derivatives of $A^a_{\ bc}$ and 
$B^a_{\ bc}$
\be
\lie A^a_{\ bc}=(\chi-\phi)A^a_{\ bc},\ \ \lie B^a_{\ bc}=(\phi-\chi)B^a_{\ bc},
\label{AyB}
\ee
 a relation which shall be used later. 

Let us now use all this information to write the Lie derivative of the connection in a convenient 
way. Note that in equation (\ref{9}) $\phi_a$ and $\chi_a$ appear projected with $P^a_{\ b}$ and 
$\Pi^a_{\ b}$ respectively suggesting that it could be interesting to write any derivative of 
$\phi$ and $\chi$ decomposed in transverse and longitudinal parts 
\be
\sphi_a\equiv \Pi_{ab}\phi^b,\ \bphi_a\equiv P_{ab}\phi^b,\ 
\schi_a\equiv P_{ab}\chi^b,\ \bchi_a\equiv\Pi_{ab}\chi^b.
\ee
If we perform this decomposition in equation (\ref{3}) and replace the terms
$\sphi_a$ and $\schi_a$ by means of (\ref{18}) we get the relation
\bea
& &\hspace{-1.5cm}\lie\left(\g^a_{\ bc}+\fr{1}{2p}(E_bP^a_{\ c}+E_cP^a_{\ b}-P_{bc}E^a)+
\fr{1}{2(n-p)}(W_b\Pi^a_{\ c}+W_c\Pi^a_{\ b}-W^a\Pi_{cb})\right)=\nonumber\\
&=&\fr{1}{2}(\bphi_bP^a_{\ c}+\bphi_cP^a_{\ b}-\bphi^aP_{bc}+\bchi_b\Pi^a_{\ c}+\bchi_c\Pi^a_{\ b}-
\bchi^a\Pi_{cb})+\fr{1}{2}(\phi-\chi)T^a_{\ bc},
\label{connection1}
\eea
but from (\ref{AyB}), (\ref{splitting-1}) and (\ref{splitting-2}) we easily deduce  
$$
\lie(A^a_{\ bc}-B^a_{\ bc})=(\chi-\phi)T^a_{\ bc},
$$
hence equation (\ref{connection1}) becomes, after some simplifications
\bea
& &\hspace{-1.5cm}
2\lie\left(\g^a_{\ bc}+\fr{1}{2p}(E_bP^a_{\ c}+E_cP^a_{\ b})+\fr{1}{2(n-p)}(W_b\Pi^a_{\ c}+W_c\Pi^a_{\ b})+
\fr{1}{2}(P^a_{\ p}-\Pi^a_{\ p})M^p_{\ cb}\right)=\nonumber\\
&=&\bphi_bP^a_{\ c}+\bphi_cP^a_{\ b}-\bphi^aP_{cb}+\bchi_b\Pi^a_{\ c}+\bchi_c\Pi^a_{\ b}-\bchi^a\Pi_{cb}.
\label{der-connection}
\eea
The geometric object inside the Lie derivative, denoted by $\bbg^a_{\ bc}$,
 is the sum of the metric connection $\g^a_{\ bc}$ plus the rank-3 tensor 
\be
\hspace{-1cm}L^a_{\ bc}\equiv\fr{1}{2p}(E_bP^a_{\ c}+E_cP^a_{\ b})+\fr{1}{2(n-p)}(W_b\Pi^a_{\ c}+W_c\Pi^a_{\ b})+
\fr{1}{2}(P^a_{\ p}-\Pi^a_{\ p})M^p_{\ bc},
\label{labc}
\ee
so it is clear that it represents a new linear connection. 
As we will see during our calculations this linear connection is fully adapted to the calculations involving 
bi-conformal vector fields and it will be extensively used in this paper. 

\begin{defn}[bi-conformal connection]
The connection whose com\-po\-ne\-nts are given by $\bbg^a_{\ bc}$ is called bi-conformal connection. 
The covariant derivative constructed from the bi-conformal connection shall be denoted by $\bnb$ and the 
curvature tensor constructed from it by $\bar{R}^a_{\ bcd}$.
\label{biconformal}
\end{defn}

Since $L^a_{\ bc}$ is symmetric in the indexes $bc$, we see that the bi-conformal connection is symmetric 
so it has no torsion and all
the identities involving only the covariant derivative $\bnb$ or the curvature $\bR^a_{\ bcd}$ remain the same
as for the case of a metric connection. However, this connection does not in general stem from a metric tensor 
as can be seen in explicit examples.
This means that certain properties of the curvature tensor of a metric 
connection are not true for $\bR^a_{\ bcd}$. 
We recall that for a symmetric connection the Riemann tensor is only antisymmetric in the last pair of indexes.

Connections defined in terms of $P_{ab}$, $\Pi_{ab}$ and their covariant derivatives as in (\ref{labc}) 
have been already considered in the literature in relation with the study of integrable 
and {\em parallel} distributions \cite{WALKER1,WALKER2,WILLMORE1,WILLMORE2} (see also \cite{REINHART}).
A distribution can be univocally fixed by means of two orthogonal and complementary projectors 
$P^a_{\ b}$, $\Pi^a_{\ b}$ and so the study of these projectors can give us information about the
geometric properties of the distribution.  

It is not very difficult to derive an identity relating the curvature tensor calculated from the 
bi-con\-for\-mal connection and the 
curvature tensor associated to the connection $\g^a_{\ bc}$
\be
\bR^a_{\ bcd}=R^a_{\ bcd}+2\nb_{[c}L^a_{\ d]b}+2L^a_{\ r[c}L^r_{\ d]b},
\label{curvatures-relation}
\ee
being this a thoroughly general identity for two symmetric connections 
$\bbg^a_{\ bc}$ and $\g^a_{\ bc}$ differing in a tensor $L^a_{\ bc}$ (\cite{SCHOUTEN}, p. 141).

The relation between the covariant derivatives $\bnb$ and $\nb$ acting on any tensor 
$X^{a_1\dots a_p}_{\ b_1\dots b_q}$ is 
\be
\hspace{-1cm}\bnb_aX^{a_1\dots a_r}_{\ b_1\dots b_q}=\nb_aX^{a_1\dots a_r}_{\ b_1\dots b_q}+\sum_{s=1}^rL^{a_s}_{\ ac}
X^{a_1\dots a_{s-1}ca_{s+1}\dots a_r}_{\ b_1\dots b_q}
-\sum_{s=1}^qL^{c}_{ab_s}X^{a_1\dots a_r}_{\ b_1\dots b_{s-1}cb_{s+1}\dots b_q},
\label{barcon-bar}
\ee
which again has general validity for two symmetric connections whose difference is a tensor $L^a_{\ bc}$ 
\cite{EISENHART}.
As a first application of this identity we may compare the covariant derivatives of the tensor $L^a_{\ bc}$ 
which leads us to the identity
$$
\bnb_{[a}L^b_{\ c]d}=\nb_{[a}L^b_{\ c]d}+2L^b_{\ r[a}L^r_{\ c]d},
$$ 
from which we can rewrite (\ref{curvatures-relation}) in terms of $\bnb$  
\be
R^a_{\ bcd}=\bR^a_{\ bcd}-2\bnb_{[c}L^a_{\ d]b}+2L^a_{\ r[c}L^r_{\ d]b}.
\label{curvature-connection}
\ee
Of course this last equation could have been obtained
 from (\ref{curvatures-relation}) by means of the replacements
$R^a_{\ bcd}\leftrightarrow\bR^a_{\ bcd}$, 
$L^a_{\ bc}\rightarrow -L^a_{\ bc}$ and $\nb_a\rightarrow\bnb_a$.
\begin{exmp}
\label{example1}
{\em
To realize the importance of the bi-conformal connection in future calculations, let us calculate 
its components for a {\em conformally separable} pseudo-Riemannian manifold 
(see definition \ref{separable}) given in local coordinates $x\equiv\{x^1,\dots,x^n\}$ by
\be
\hspace{-.6cm}
ds^2=\Xi_1(x^1,\dots,x^n)G_{\alpha\beta}(x^{\delta})dx^{\alpha}dx^{\beta}
+\Xi_2(x^1,\dots,x^n)G_{AB}(x^C)dx^Adx^B.
\ee
Here Greek indexes range from 1 to $p$ and uppercase Latin indexes from $p+1$ to $n$ so the 
metric tensors $G_{\alpha\beta}$ and $G_{AB}$ are of rank $p$ and $n-p$ respectively.
The tensors $G^{\alpha\beta}$ and $G^{AB}$ are defined in the obvious way and they are 
used to raise Greek and uppercase Latin indexes respectively.
The non-zero Christoffel symbols for this metric are
\bnr
\Chr^{\alpha}_{\ \beta\g}=\fr{1}{2\Xi_1}G^{\alpha\rho}(\d_{\beta}(\Xi_1G_{\alpha\rho})+
\d_{\g}(\Xi_1G_{\rho\beta})-\d_{\rho}(\Xi_1G_{\beta\g})),\\
\Chr^{A}_{\ BC}=\fr{1}{2\Xi_1}G^{AD}(\d_{B}(\Xi_1G_{CD})+
\d_{C}(\Xi_1G_{DB})-\d_{D}(\Xi_1G_{BC})),\\
\Chr^{\alpha}_{\beta A}=\fr{1}{2\Xi_1}\delta^{\alpha}_{\ \beta}\d_A\Xi_1,\ \ 
\Chr^A_{\ B\alpha}=\fr{1}{2\Xi_2}\delta^A_{\ B}\d_{\alpha}\Xi_2,
\enr
from which the only nonvanishing components of $M_{abc}$, $E_a$, $W_a$ are
\bea
M_{\alpha AB}=\d_{\alpha}(\Xi_2G_{AB}),\ \ 
M_{A\alpha\beta}=-\d_{A}(\Xi_1G_{\alpha\beta}),\nonumber\\
E_{A}=-\d_{A}\log|\mbox{det}(\Xi_1G_{\alpha\beta})|,\ \ 
W_{\alpha}=-\d_{\alpha}\log|\mbox{det}(\Xi_2G_{AB})|.
\label{27}
\eea
Therefore we get for the components of the bi-conformal connection
\bea
\bar{\Chr}^{\alpha}_{\ \beta\phi}&=&\fr{1}{2\Xi_1}(\delta^{\alpha}_{\ \beta}\d_{\phi}\Xi_1+\delta^{\alpha}_{\ \phi}\d_{\beta}\Xi_1-
G^{\alpha\rho}G_{\beta\phi}\d_{\rho}\Xi_1)+\Chr^{\alpha}_{\ \beta\phi}(G),\nonumber\\
\bar{\Chr}^A_{\ BC}&=&\fr{1}{2\Xi_2}(\delta^{A}_{\ B}\d_{C}\Xi_2+\delta^{A}_{\ C}\d_{B}\Xi_2-
G^{AR}G_{BC}\d_{R}\Xi_2)+\Chr^{A}_{\ BC}(G)\nonumber\\
\bar{\Chr}^{\alpha}_{\ \beta C}&=&\bar{\Chr}^{A}_{\ B\phi}=0,
\label{new-components}
\eea
where $\Chr^{\alpha}_{\ \beta\phi}(G)$ and $\Chr^{A}_{\ BC}(G)$ are the Christoffel symbols of the metrics 
$G_{\alpha\beta}$ and $G_{AB}$ respectively. From the above formulae we deduce that the bi-conformal connection is 
fully adapted to a conformally separable pseudo-Riemannian manifold because its components clearly split in two 
parts being each of them the Christoffel symbols of the metrics $G_{\alpha\beta}$, $G_{AB}$ plus terms 
involving the derivatives of the factors $\Xi_1$ and $\Xi_2$. We will take advantage of this property in 
section \ref{sufficient}
where we will find  an invariant local characterization of conformally separable 
pseudo-Riemannian manifolds.} 
\end{exmp}

We calculate next the covariant derivative with respect to the bi-conformal connection of a number of 
tensors.

\begin{prop}
The following identities hold true 
\bea
\bnb_aP_{bc}&=&\nb_aP_{bc}-\fr{1}{p}E_aP_{bc}-\fr{1}{2p}(E_bP_{ac}+E_cP_{ab})-\nonumber\\ 
&-&\fr{1}{2}(P_{cp}M^p_{\ ab}+P_{bp}M^p_{\ ac}),\label{identity-1}\\
2\bnb_aP^b_{\ c}&=&2\nb_aP^b_{\ c}+P^{bq}P^r_{\ c}M_{qra}-\Pi^{bq}P^r_{\ c}M_{qra}-
P^b_{\ q}M^q_{\ ac}+\fr{1}{n-p}W_c\Pi^b_{\ a}-\nonumber\\
&-&\fr{1}{p}E_cP^b_{\ a},
\label{identity-2}\\
\bnb_aP^{bc}&=&\nb_aP^{bc}+\fr{1}{p}E_aP^{bc}+\fr{1}{2(n-p)}(W^c\Pi^b_{\ a}+W^b\Pi^c_{\ a})-\nonumber\\
&-&\fr{1}{2}(M^b_{\ ar}P^{rc}+M^c_{\ ar}P^{rb}),
\label{identity-3}
\eea
and all the identities formed with the replacements $P_{ab}\rightarrow\Pi_{ab}$, $p\rightarrow n-p$.
\label{identity-set}
\end{prop}
\begin{pf}
All these identities are proven by means of (\ref{barcon-bar}) and the use of  
 properties (\ref{vwe}).\qed 
\end{pf}
Using the above properties we can get more interesting identities to be used later on.
\bea
\bnb_aP^{ab}=\bnb_a\Pi^{ab}=\bnb_aP^a_{\ b}=\bnb_a\Pi^a_{\ b}=0,\label{contraction}\\
P^{bc}\bnb_aP_{bc}=-E_a,\ \ P_{bc}\bnb_aP^{bc}=E_a,\label{EW}\\
\Pi^{bc}\bnb_a\Pi_{bc}=-W_a,\ \ \Pi_{bc}\bnb_a\Pi^{bc}=W_a,\label{WE}\\
P^d_{\ r}\bnb_b\Pi^r_{\ d}=\Pi^d_{\ r}\bnb_bP^r_{\ d}=P^{dr}\bnb_{d}P_{rb}=\Pi^{dr}\bnb_d\Pi_{rb}=0.\label{PPi}
\eea
Note that index raising and lowering do not commute with $\bnb$ so we must be very careful when we 
raise or lower indexes in tensor expressions involving $\bnb$.

\section{Normal form and dimension of maximal Lie algebras of bi-con\-for\-mal fields}
\label{normal-dimension}
We turn now our attention to the calculation of the full normal form coming from 
the differential conditions (\ref{bi-conformal}). A detailed explanation of the general procedure and 
relevance of this calculation for a general symmetry can be found in \cite{EISENHARTII,BI-CONFORMAL} 
(see also \cite{YANO} 
for the calculation in the cases of the most studied symmetries in General Relativity such as isometries 
and conformal motions). Before starting the 
calculation and for the sake of completeness let us give a very brief sketch of the whole procedure. 
We must differentiate the condition (\ref{bi-conformal}) a number of times in such a way that we get enough 
equations to isolate the derivatives of certain variables (system variables) in terms of themselves (this is 
achieved typically by means of the resolution of a linear system of equations). 
The so obtained derivatives give rise to the normal form associated to our symmetry. 
Along the differentiating process one
may obtain equations whose linear combinations no longer contain derivatives of the system variables 
(constraints). Examples of such constraints in our case are the differential conditions themselves and 
(\ref{18}) (actually these are the only constraints as we will show in $\S$\ref{constraints}). 
   
It is possible to meet cases in which the normal form cannot be achieved. This means that 
one cannot obtain enough equations to isolate all the derivatives obtained through the derivation process.
The main implication of this is that the Lie algebra of vector fields fulfilling the starting differential 
condition is infinite dimensional as opposed to the case in which there is such a normal form. Therefore the 
calculation of the normal form allows us to tell apart the cases with an infinite dimensional Lie algebra of 
vector fields from those representing finite dimensional Lie algebras. In the latter case we can even go
further and determine the highest dimension of these Lie algebras as the total number of system variables 
minus the number of linearly independent constraints.  

We start out our calculation with the substitution of (\ref{der-connection}) into (\ref{lie-xi}) which yields
\bea
\bnb_b\Psi_c^{\ a}&+&\xi^d\bar{R}^a_{\ cdb}=\fr{1}{2}(\bphi_bP^a_{\ c}+\bphi_cP^a_{\ b}-\bphi^aP_{cb}+
\bchi_b\Pi^a_{\ c}+\bchi_c\Pi^a_{\ b}-\bchi^a\Pi_{cb}),\nonumber\\
\Psi_c^{\ a}&\equiv&\bnb_c\xi^a.
\label{psi}
\eea
Next we replace in (\ref{lie-curvature}) the Lie derivatives of the bi-conformal connection by their 
expressions given by (\ref{der-connection}) getting 
\bea
\lie\bR^d_{\ cab}=\bnb_{[a}\bphi_{b]}P^d_{\ c}&+&P^d_{\ [b}\bnb_{a]}\bphi_c-P_{c[b}\bnb_{a]}\bphi^d+
\bnb_{[a}\bchi_{b]}\Pi^d_{\ c}+\Pi^d_{\ [b}\bnb_{a]}\bchi_c-\nonumber\\
-\Pi_{c[b}\bnb_{a]}\bchi^d+\bphi_{[b}\bnb_{a]}P^d_{\ c}&+&\bphi_c\bnb_{[a}P^d_{\ b]}-\bphi^d\bnb_{[a}P_{b]c}+
\bchi_{[b}\bnb_{a]}\Pi^d_{\ c}+\bchi_c\bnb_{[a}\Pi^d_{\ b]}-\nonumber\\
-\bchi^d\bnb_{[a}\Pi_{b]c}\label{lie-curvatura}.
\eea
The game is now to isolate from this expression $\bnb_a\bphi_b$ and $\bnb_{a}\bchi_b$ (these rank-2 tensors are 
not symmetric in general). 
The forthcoming calculations split in two groups which are dual under the interchange
$P_{ab}\Leftrightarrow\Pi_{ab}$, $p\Leftrightarrow n-p$ (only the calculations of the first group are shown).
Multiplying (\ref{lie-curvatura}) by $P^a_{\ r}$ we obtain
\bea
& &\lie(P^d_{\ r}\bR^r_{\ cab})=P^d_{\ c}\bnb_{[a}\bphi_{b]}+P^d_{\ [b}\bnb_{a]}\bphi_c-
P^d_{\ r}P_{c[b}\bnb_{a]}\bphi^r-P^d_{\ r}\Pi_{c[b}\bnb_{a]}\bchi^r+\nonumber\\
&+&P^d_{\ r}\bphi_{[b}\bnb_{a]}P^r_{\ c}+\bphi_cP^d_{\ r}\bnb_{[a}P^r_{\ b]}-
\bphi^d\bnb_{[a}P_{b]c}+P^d_{\ r}\bchi_{[b}\bnb_{a]}\Pi^r_{\ c}+\bchi_cP^d_{\ r}\bnb_{[a}\Pi^r_{\ b]}.\nonumber\\
\label{seis} 
\eea
Contraction of the indexes $d$-$c$ in the above expression yields
\be
\bnb_a\bphi_b=\bnb_b\bphi_a+\fr{2}{p}\lie(P^d_{\ r}\bR^r_{\ dab}),
\label{nueve}
\ee 
while the contraction of indexes $d$-$a$ and use of identities (\ref{contraction})-(\ref{PPi}) entails
\be
2\lie(P^d_{\ r}\bR^r_{\ cdb})=P^d_{\ c}\bnb_d\bphi_b+P^d_b\bnb_d\bphi_c-\bphi^r\bnb_rP_{bc}-\bnb_a\bphi^aP_{bc}-
p\bnb_b\bphi_c.
\label{ocho}
\ee
Equations (\ref{nueve}) and (\ref{ocho}) can now be combined in a single expression which is 
\bea
2\lie\left[P^d_{\ r}\bR^r_{\ cdb}-\fr{1}{p}(P^d_{\ c}P^r_{\ q}\bR^q_{\ rdb}+P^d_{\ b}P^r_{\ q}\bR^q_{\ rdc}-
P^r_{\ q}\bR^q_{\ rbc})\right]=\nonumber\\
=(2-p)\bnb_b\bphi_c-\bnb_a\bphi^aP_{bc}-\bphi_d(\bnb_bP^d_{\ c}+\bnb_cP^d_{\ b}+P^{rd}\bnb_rP_{bc}).
\label{diez}
\eea
Multiplying here with $P^{cb}$ leads, after a little bit of algebra, to  
\be
\bnb_a\bphi^a=\fr{1}{1-p}(\lie\bR^0+\phi\bR^0),\ \ p\neq 1,\ \ \bR^0=P^d_{\ r}\bR^r_{\ cdb}P^{cb},
\label{once}
\ee

On the other hand (\ref{diez}) can be further simplified if we take into account the identity
\bnr
\bnb_bP^d_{\ c}+\bnb_cP^d_{\ b}-P^{rd}\bnb_rP_{bc}
&=&\Pi^{rd}\nb_rP_{bc}+\fr{1}{2}(\Pi^{dq}\Pi^r_{\ c}M_{qrb}+\Pi^{dq}\Pi^r_{\ b}M_{qrc})+\\
&+&\fr{1}{2(n-p)}(W_c\Pi^d_{\ b}+W_b\Pi^d_{\ c}),
\enr
which is easily derived by writing all the covariant derivatives with respect to the bi-conformal connection 
of the projectors in terms of ordinary covariant derivatives (proposition \ref{identity-set}). 
So plugging (\ref{once}) into (\ref{diez}) we get
\be
(2-p)\bnb_b\bphi_c=\lie L^0_{bc}+2\bphi^r\bnb_rP_{bc},
\label{normal-1}
\ee
where
\be
\hspace{-1cm}L^0_{\ bc}\equiv 
2\left[P^d_{\ r}\bR^r_{\ cdb}-\fr{1}{p}(P^d_{\ c}P^r_{\ q}\bR^q_{\ rdb}+P^d_{\ b}P^r_{\ q}\bR^q_{\ rdc}-
P^r_{\ q}\bR^q_{\ rbc})\right]+\fr{\bR^0}{1-p}P_{bc}.
\label{L0}
\ee  
The duals of (\ref{normal-1}), (\ref{L0}) are 
\be
(2-n+p)\bnb_b\bchi_c=\lie L^1_{bc}+2\bchi^r\bnb_r\Pi_{bc},
\label{normal-2}
\ee
and 
\bea
L^1_{\ bc}\equiv 
2\left[\Pi^d_{\ r}\bR^r_{\ cdb}-\fr{1}{n-p}(\Pi^d_{\ c}\Pi^r_{\ q}\bR^q_{\ rdb}+\Pi^d_{\ b}\Pi^r_{\ q}\bR^q_{\ rdc}-
\Pi^r_{\ q}\bR^q_{\ rbc})\right]&+&\nonumber\\
+\fr{\bR^1}{1-n+p}\Pi_{bc},\ \ \bR^1\equiv\Pi^d_{\ r}\bR^r_{\ cdb}\Pi^{cb}.& &
\label{L1}
\eea
To complete the normal form we need now the derivatives of $\sphi_a$ and $\schi_a$ which are obtained through 
the differentiation of (\ref{18}) (identity (\ref{lie-conmmutation}) must be used to get these derivatives)
\bea
-p\bnb_b\sphi_a=\lie(\bnb_bE_a)+\fr{1}{2}(\bchi_bE_a+\bchi_aE_b-(\bchi^rE_r)\Pi_{ab})\\
(p-n)\bnb_b\schi_a=\lie(\bnb_bW_a)+\fr{1}{2}(\bphi_bW_a+\bphi_aW_b-(\bphi^rW_r)P_{ab})
\eea

\subsection{Normal form of the differential conditions} 
The above calculations give us the sought normal form for the differential conditions 
(\ref{bi-conformal}) being these gathered in the  following set of equations  
\bea
(a)\ \ \bnb_a\phi&=&\bphi_a+\sphi_a,\ \bnb_a\chi=\bchi_a+\schi_a,\nonumber\\
(b)\ \ \bnb_b\sphi_a&=&\fr{-1}{p}\left[\lie(\bnb_bE_a)+\fr{1}{2}(\bchi_bE_a+\bchi_aE_b-(\bchi^rE_r)\Pi_{ab})\right],\nonumber\\
(c)\ \ \bnb_b\schi_a&=&\fr{1}{p-n}\left[\lie(\bnb_bW_a)+\fr{1}{2}(\bphi_bW_a+\bphi_aW_b-(\bphi^rW_r)P_{ab})\right],\nonumber\\
(d)\ \ \bnb_b\bphi_c&=&\fr{1}{2-p}\left[\lie L^0_{bc}+2\bphi^r\bnb_rP_{bc}\right],\label{normal-form}\\
(e)\ \ \bnb_b\bchi_c&=&\fr{1}{2-n+p}\left[\lie L^1_{bc}+2\bchi^r\bnb_r\Pi_{bc}\right],\nonumber\\
(f)\ \ \bnb_b\xi^a&=&\Psi_b^{\ a},\nonumber\\ 
(g)\ \ \bnb_b\Psi_c^{\ a}&=&\fr{1}{2}(\bphi_bP^a_{\ c}+\bphi_cP^a_{\ b}-\bphi^aP_{cb}+
\bchi_b\Pi^a_{\ c}+\bchi_c\Pi^a_{\ b}-\bchi^a\Pi_{cb})-\xi^d\bar{R}^a_{\ cdb}.\nonumber 
\eea

A first glance at these equations reveals us that this normal form does not always exist. 
To be precise if either $p=2$ or $p=n-2$ the derivatives $\bnb_a\bphi_b$ and $\bnb_a\bchi_b$ cannot be isolated
and the system cannot be ``closed'' (in fact these derivatives cannot be isolated even if we perform further 
derivatives of any of the above equations). We must also remember at this point that the tensors $L^0_{ab}$ and 
$L^1_{ab}$ are well-defined unless $p=1$ and $p=n-1$ respectively (see equations (\ref{L0}) and (\ref{L1})).
Therefore we have proven the following theorem (compare to proposition 6.2 of \cite{BI-CONFORMAL})
\begin{thm}
The only cases in which the Lie algebra $\Gl(S)$ can be infinite dimensional occur if 
and only if $p=1$, $p=2$, $p=n-1$, $p=n-2$.
\end{thm}
The result of this theorem is intuitively clear if we realize that bi-conformal vector fields are somehow
conformal motions for the projectors $P_{ab}$ and $\Pi_{ab}$. Therefore if any of them projects onto one or 
two dimensional vector spaces the associated Lie algebras may turn out to be infinite dimensional as we have just found.   

Equations (\ref{normal-form}) are no longer a normal form system 
if some of the derivatives involved 
vanish. This happens for instance if part of the gauge functions are constants or their second covariant
derivatives with respect to the bi-conformal connection   
are zero. In all this work we will assume that these derivatives are not zero 
in an open neighbourhood of a point
 leaving the study of any other cases for a forthcoming publication. 

The normal form can be used to establish the minimum 
conditions under which bi-conformal
vector fields are smooth vector fields.

\begin{prop}
Let $\xiv$ be a bi-conformal vector field at least $C^2$ in a neighbourhood $\U_x$ of a point $x$ belonging to a 
manifold $V$ with a $C^{\infty}$ metric tensor. 
If $\phi$, $\chi$ are at least $C^2$ on $\U_x$ then $\xiv\in C^{\infty}(\U_x)$. 
\label{differentiability}
\end{prop}  
\begin{pf} To prove this result it is enough to show that the covariant derivatives of $\xiv$ with res\-pect to 
the bi-conformal connection exist at any order. The first and second derivatives of $\xiv$ are equations 
(\ref{normal-form})-{\em f} and (\ref{normal-form})-{\em g} and higher derivatives are calculated from 
this last equation. When we derive (\ref{normal-form})-{\em g} we need 
$\bnb_b\bphi_c$ and $\bnb_b\bchi_c$ which exist on $\U_x$ as $\phi$, $\chi\in C^2(\U_x)$ 
and these are obtained through
equations (\ref{normal-form})-{\em d} and (\ref{normal-form})-{\em e} in which only
derivatives of $\xiv$, $\bphi_a$ and $\bchi_a$ of order less or equal than one appear.
This makes clear that no other equation of (\ref{normal-form}) but the ones mentioned so far are involved 
in the calculation of the derivatives of $\xiv$ and so we can obtain them in as high order as we wish.
 \qed
\end{pf}

Related to this is the following result (the technique used in the proof
has been employed in \cite{LIBROHALL} for other symmetries in the 
framework of General Relativity).
\begin{prop}
Under the hypotheses of previous proposition if a bi-con\-for\-mal vector field $\xiv$ is such that 
$\xi^a|_x=0$, $\bnb_b\xi^a|_x=0$, $\bnb_c\bnb_b\xi^a|_x=0$ then $\xiv\equiv 0$ in a neighbourhood of $x$.
\end{prop}
\begin{pf}
Evaluation of the last equation of (\ref{normal-form}) at $x$ entails
$$
(\bphi_bP^a_{\ c}+\bphi_cP^a_{\ b}-\bphi^aP_{cb}+
\bchi_b\Pi^a_{\ c}+\bchi_c\Pi^a_{\ b}-\bchi^a\Pi_{cb})|_x=0,
$$
from which, projecting down with $P^{ab}$ and $\Pi^{ab}$, we deduce that $\bphi_a|_x=\bchi_a|_x=0$.
Now let $\g(t)$ be a smooth curve on $V$ such that $\g(0)=x$ with $\g(t)$ lying in
a coordinate neighbourhood of $x$ for all $t$ in the interval $(-\epsilon,\epsilon)$. If we denote by $\dot{\g}^a(t)$ the tangent 
vector to this curve we may define the derivatives
$$
\fr{\bar{D}\bphi_a}{dt}\equiv\dot{\g}^r\bnb_r\bphi_a,\ \fr{\bar{D}\bchi_a}{dt}\equiv\dot{\g}^r\bnb_r\bchi_a,\ 
\fr{\bar{D}\xi^a}{dt}\equiv\dot{\g}^r\bnb_r\xi^a,\ \fr{\bar{D}\Psi_c^{\ a}}{dt}\equiv\dot{\g}^r\bnb_r\Psi_c^{\ a}, 
$$   
where all quantities are evaluated on $\g(t)$. Contracting (\ref{normal-form})-{\em d}-(\ref{normal-form})-{\em g} with 
$\dot{\g}^b$ we can transform these equations into a first order ODE system in the variables $\bphi_a(\g(t))$, $\bchi_a(\g(t))$,
$\xi^a(\g(t))$ and $\Psi_a^{\ b}(\g(t))$ with $\bar{D}/dt$ as derivation. 
From the above all these variables vanish at $t=0$ so according to 
the standard theorem of uniqueness for ODE systems the variables are identically zero along the curve $\g(t)$ (and in particular 
$\xiv|_{\g(t)}=0$). As $\g(t)$ was chosen arbitrarily we conclude that $\xiv\equiv 0$ in a whole neighbourhood of $x$.\qed 
\end{pf}
\begin{rem}{\em
If the manifold is connected (and hence path connected) then 
the curve $\g(t)$ can be chosen joining any pair of points of the manifold
and then the vector field $\xiv$ is zero everywhere and not just in 
a single neighbourhood of a point.} 
\end{rem}

From the calculations performed in this proof and in the proof of proposition \ref{differentiability}  
we deduce as a simple corollary that 
$\bnb^{(m)}\xi^a|_x\equiv 0$ for all $m\in{\mathbb N}$ if it holds for $m=1,2$.  

\subsection{Constraints}
\label{constraints}
From the above calculation of the normal form, the system variables are read off at once. These are the variables 
appearing under derivation in the left hand side of (\ref{normal-form}). However, they are not algebraically independent
because in the calculation process of (\ref{normal-form}) some of the equations involved do not contain derivatives 
of the variables at all (system constraints). The most evident case of these constraints are the differential 
conditions (\ref{bi-conformal}) themselves. If we review the whole procedure followed to get (\ref{normal-form}) we 
deduce that the other set of constraints between the system variables is (\ref{18}) so (\ref{normal-form}) must be 
complemented with 
\be
(I)\left.\begin{array}{c}
\lie P_{ab}=\phi P_{ab}\\
\lie\Pi_{ab}=\chi\Pi_{ab}
\end{array}\right\},\ \ 
(II)\left.\begin{array}{l}
\lie E_a=-p\sphi_a\\
\lie W_a=-(n-p)\schi_a\end{array}\right\}
\label{ligaduras}
\ee 
In order to clarify that these two sets of equations are truly the constraints associated with (\ref{normal-form}) 
we must show that they arise as a linear combination of some of the higher covariant derivatives  of 
(\ref{bi-conformal}) employed to get the normal form. 
An equation equivalent to the first derivative of 
(\ref{bi-conformal}) is (\ref{9}) but only its projections by 
$P^{ab}$ and $\Pi^{ab}$ (equation \ref{18}) really matter to work out the normal form and these 
are (\ref{ligaduras})-II. Equation (\ref{psi}) is also obtained from (\ref{9}) and it is part of the normal 
system and not a constraint.
 As for the other derivatives they do not give rise to 
any more equations with no derivatives of the system variables so the above equations are the only constraints 
we must care about.   

A first application of all the above calculations comes in the following result, 
already proven in \cite{BI-CONFORMAL}  
using a normal form system written in terms of different variables.
\begin{thm}
If the Lie algebra $\Gl(S)$ is finite dimensional then its dimension is bo\-un\-ded from above by 
$N=p(p+1)/2+(n-p)(n-p+1)/2$.
\label{bounded}
\end{thm}

\begin{pf}
To prove this theorem we must state what the upper bound to the maximum number of {\em integration constants} for the system 
(\ref{normal-form}) is. As is very well known from the theory of normal systems of PDE's (see e. g. \cite{EISENHARTII}) 
such number is the number of system variables minus the number of linearly independent constraints. The following 
table summarizes these numbers for our system.

\begin{table}[h]
\caption{Calculation of the highest dimension of $\Gl(S)$}
\label{tab:1}       
\centering
\begin{tabular}{ccccc|cc}
\hline\hline\noalign{\smallskip}
\multicolumn{5}{c}{system variables} & 
\multicolumn{2}{c}{Constraints}\\
\noalign{\smallskip}\hline\hline\noalign{\smallskip}
$\phi,\chi$ & $\sphi_a$, $\bphi_a$ & $\schi_a$, $\bchi_a$ & $\xi^a$ & $\Psi_a^{\ b}$ & 
eq. (\ref{ligaduras})-I & eq. (\ref{ligaduras})-II\\
\noalign{\smallskip}\hline\noalign{\smallskip}
 2 & $n$ & $n$ & $n$ & $n^2$ & $n(n+1)/2+p(n-p)$ & $n$\\
\end{tabular}
\end{table}

We have written explicitly the system variables and the total number for each of them. The constraints are also 
indicated together with how many linearly independent equations each constraint amounts to. This last part is not 
evident as opposed to the co\-un\-ting for the system variables so the rest of the proof is devoted to show 
that the numbers given in table \ref{tab:1} for the constraint
equations are indeed correct.  

\medskip
\noindent
{\em Equations (\ref{ligaduras})-I}. First of all, we expand the Lie derivatives of these equations
\bea 
\xi^c\bnb_cP_{ab}+\Psi_p^{\ c}(\dd^p_{\ a}P_{bc}+\dd^p_{\ b}P_{ac})=\phi P_{ab},\nonumber\\ 
\xi^c\bnb_c\Pi_{ab}+\Psi_p^{\ c}(\dd^p_{\ a}\Pi_{bc}+\dd^p_{\ b}\Pi_{ac})=\chi\Pi_{ab},
\label{expansion}
\eea
where the standard definition of the Lie derivative of a tensor $P_{ab}$ has been applied
\be
\lie P_{ab}\equiv\xi^c\bnb_c P_{ab}+\bnb_a\xi^cP_{cb}+\bnb_b\xi^cP_{ac},
\ee   
(observe that the general formula of the Lie derivative with respect to a vector in terms of its covariant
derivatives still holds under a symmetric connection). We define new indexes 
$A$, $B$, $B'$ in such a way that
$$
P_A\equiv P_{a b},\ \Psi_{B}\equiv\bnb_p\xi^q,\ \xi^{B'}=\xi^c, 
$$
so capital indexes group together certain combinations of small indexes (explicitly $A=\{a,b\}$, $B=\{p, q\}$ and
$B'=c$. The ranges of the new indexes are 
$A=1,\dots, n(n+1)/2$, $B=1,\dots n^2$, $B'=1,\dots n$.
Using these new labels we can write in
matrix notation the homogeneous system posed by these constraints (we only concentrate in 
the first of (\ref{expansion}))
\be
(\begin{array}{ccc}
   M_A^{\ B} & (\bnb P)_{B'A} & P_A 
\end{array}) \left(\begin{array}{c}
                              \Psi_B\\
                              \xi^{B'}\\
                              -\phi \end{array}\right)=0, 
\label{matrix-form}
\ee
($A$= row index, $B$, $B'$=column indexes) 
where the explicit expressions of the matrices read
$$
(\bnb P)_{B'A}=\bnb_cP_{ab},\ 
\ M_A^{\ B}=\dd^p_{\ a}P_{bq}+\dd^p_{\ b}P_{aq}.
$$

The number of linearly independent equations is just the rank of the 
matrix system of (\ref{matrix-form}). 
In principle the rank of this matrix will depend on the projector $P_{ab}$ and its covariant 
derivative, 
meaning this that it will depend on the geometry of the manifold. However, since we are interested in spaces with a 
maximum 
number of bi-conformal vector fields it is enough to find the least rank of the above matrix for all possible 
projectors $P_{ab}$. We start first studying the rank of $M_A^{\ B}$ whose  
 nonvanishing components occur in the following cases (no summation 
over the repeated indexes)
\begin{eqnarray*}
p=a,\ q=b &\Rightarrow & \left\{\begin{array}{cc}
M_A^{\ A}=\dd^a_{\ a}P_{bb}, & a\neq b,\\
M_A^{\ A}=2 P_{aa},\ a=b. & 
       \end{array}\right.\\
p=b,\ q=a &\Rightarrow &\
 M_A^{\ Q'}=\dd^b_{\ b}P_{aa},\ a\neq b,
\end{eqnarray*}
where we have assumed that we are working in the common (orthonormal) basis of eigenvectors of $P^a_{\ b}$ and 
$\Pi^a_{\ b}$ so 
$$
P^{a}_{\ b}=\mbox{diag}(\overbrace{1\dots 1}^{p} 0 \dots 0),\  
\Pi^a_{\ b}=\mbox{diag}(0\dots 0\overbrace{1\dots 1}^{n-p}).
$$ 
Hence we only need to count how many components of the type $M_{A}^{\ A}$ are different from zero 
because by construction these elements give rise to linearly independent rows of the matrix $M_A^{\ B}$
(the elements $M_{A}^{\ Q'}$ are in the same row of the matrix $M_A^{\ B}$ and they do not increase its rank because $Q'>A$).
The sought number can be obtained from the following diagram
gathering into blocks the $A$ indexes of the rows containing non-zero
elements (we express each index in terms of tensor indexes following the notation $A=(a,b)$) 
$$
\begin{array}{l|l}
\mbox{Block 1}=\{\overbrace{(1,1)\ (1,2)\ \dots\ (1,p)}^{p}\} &\ \ \mbox{Block $p+1$}=\{\overbrace{(p+1,1)\ \dots (p+1,p)}^{p}\}\\
\mbox{Block 2}=\{\overbrace{(2,2)\ (2,3)\ \dots\ (2,p)}^{p-1}\}&\ \ \mbox{Block $p+2$}=\{\overbrace{(p+2,1)\ \dots\ (p+2,p)}^p\}\\  
\mbox{Block 3}=\{\overbrace{(3,3)\ (3,4)\ \dots\ (3,p)}^{p-2}\}&\ \  \mbox{Block $p+3$}=\{\overbrace{(p+3,1)\ \dots\ (p+3,p)}^p\}\\   
\dots\dots\dots\dots\dots\dots\dots\dots\dots\dots\dots&\ \ \dots\dots\dots\dots\dots\dots\dots\dots\dots\dots\dots\dots\dots \\
\mbox{Block $p$}=\{\overbrace{(p,p)}^1\} &\ \  \mbox{Block $p+n$}=\{\overbrace{(n,1)\ \dots\ (n,p)}^p\},
\end{array}
$$
from which the rank of $M_A^{\ B}$ is 
$$
1+\dots+p+p(n-p)=\fr{1}{2}p(p+1)+p(n-p).
$$
Note that this rank only depends on algebraic properties of the projector $P_{ab}$ and not on its actual form at some 
concrete space. Addition of the matrices $(\bnb P)_{B'A}$ and $P_A$ only increase the rank of the matrix of the homogeneous system 
(\ref{matrix-form}) and so we do not need to take them into account.  The total 
number of constraints posed by (\ref{ligaduras})-I is then the rank of $M_A^{\ B}$ plus the rank of the matrix $N_A^{\ B}$ constructed
replacing $P_{ab}$ by $\Pi_{ab}$
\bnr
rank(M)+rank(N)=& & \\
=\fr{1}{2}p(p+1)+p(n-p)&+&\fr{1}{2}(n-p)(n-p+1)+p(n-p)=\nonumber\\
=\fr{1}{2}n(n+1)&+&p(n-p).
\enr 

\medskip
\noindent
{\em Equations (\ref{ligaduras})-II.} In order to perform the analysis of these 
constraints it is enough to realize that 
the 1-forms $\sphi_a$ and $\schi_a$ appearing in the right hand side of each 
equation are invariant under the projectors 
$\Pi^a_{\ b}$ and $P^a_{\ b}$ respectively. Therefore the first equation of 
(\ref{ligaduras})-II contains at least $n-p$ linearly independent 
equations and the second one $p$ being $n$ the total sum of them.

The upper bound $N$ is then 
\bnr
N=2+n+n+n+n^2-(n+\fr{1}{2}n(n+1)+p(n-p))&=&\\
=\fr{1}{2}(p+1)(p+2)+\fr{1}{2}(n-p+1)(n-p+2)& &.
\enr
\qed
\end{pf}

This proof does not guarantee the existence of a Lie algebra $\Gl(S)$ in which the dimension $N$ is attained. 
However, 
it is not difficult to find explicit examples of pseudo-Riemannian manifold possessing $N$ 
bi-conformal vector fields 

\begin{prop}
The number $N$ is the maximum dimension of $\Gl(S)$ if $p,\ n-p\notin\{1,2\}$ being
this dimension attained for any pseudo-Riemannian manifold whose line element is in local coordinates $\{x^a\}$, $a=1,\dots, n$
\be
ds^2=\phi_1^2(x^a)\eta^0_{\alpha\beta}dx^{\alpha}dx^{\beta}+\phi_2^2(x^a)\eta^1_{AB}dx^Adx^B,
\label{twopieces}
\ee
where $x^{\alpha}=\{x^1,\dots,x^p\}$, $x^A=\{x^{p+1},\dots,x^{n}\}$ are sets of coordinates and
$\bfeta^0$, $\bfeta^1$ flat metrics of the appropriate signatures
depending only on the coordinates $\{x^{\alpha}\}$ and $\{x^{A}\}$ respectively.
\label{canonical-form}
\end{prop}

\begin{pf} 
This result is proposition 6.1 of \cite{BI-CONFORMAL}. \qed
\end{pf}
Spaces of previous proposition are called {\em bi-conformally flat}. As stated in definition 
\ref{classification} they are a particular case of 
{\em conformally separable spaces} and we may ask if they are the only pseudo-Riemannian 
manifolds admitting $N$ independent bi-conformal vector fields. The answer to this and other 
questions such as their geometric characterization can be settled by calculating the complete
 integrability conditions of (\ref{normal-form}). Remarkably this has been already done but 
the whole procedure 
relies on hefty algebraic manipulations so we have preferred to present these results in a 
subsequent paper. The interested reader can find the full details of these calculations in \cite{GRQC}.

\section{Local geometric characterization of conformally separable pseudo-
Rie\-man\-ni\-an manifolds}
\label{sufficient}
In this section we will show how the bi-conformal connection can be used to 
derive an invariant geometric characterization of conformally separable pseudo-Riemannian
manifolds. To begin with we define in precise terms what a conformally separable pseudo-Riemannian manifold is 
\begin{defn}
The pseudo-Riemannian manifold $(V,\G)$ is said to be conformally separable at the point $q\in V$ if there exists a 
local coordinate chart 
$x\equiv\{x^1,\dots,x^n\}$ based at $q$ in which the metric tensor takes the form
\be
\rmg_{ab}(x)=\left\{
\begin{array}{c}
\Xi_1(x)G_{\alpha\beta}(x^\g),\ 1\leq\alpha,\beta,\g\leq p\\
\Xi_2(x)G_{AB}(x^C),\ p+1\leq A,B,C\leq n\\
0\ \ \mbox{otherwise.}
\end{array}\right.
\label{metric-separable}
\ee
where $\Xi_1$, $\Xi_2$ are $C^2$ functions on the open set defining the coordinate chart.
$(V,\G)$ is conformally separable if it is so at every point $p\in V$. Any of the metric tensors 
$\Xi_1G_{\alpha\beta}$, $\Xi_2G_{AB}$ shall be called a
leaf metric.
\label{separable}
\end{defn}
Henceforth all our results will deal with conformally separable pseudo-Rieman\-nian manifolds 
at a point. From now on 
when working with conformally separable spaces written in the form of (\ref{metric-separable})
we adopt the convention that 
Greek letters label indexes associated to one of the leaf metrics whereas uppercase Latin characters 
are used for the other one. 

Conformally separable pseudo-Riemannian manifolds are also known as 
{\em double twisted products}. They comprise a number of particular
cases which have received wide attention in the literature under different
nomenclatures. A summary of them is presented next. 
\begin{defn} 
Let $x\equiv\{x^a\}$, $a=1\dots n$ be the local coordinate system introduced
in definition \ref{separable}.  
A conformally separable manifold can then be classified in terms of the form
that the functions $\Xi_1$, $\Xi_2$ take in these coordinates as 
\begin{enumerate}
\item decomposable or reducible:\\
 $ds^2=G_{\alpha\beta}(x^{\epsilon})dx^{\alpha}dx^{\beta}+G_{AB}(x^C)dx^Adx^B$, 
\item semi-decomposable, semi-reducible or warped product:\\ 
$ds^2=G_{\alpha\beta}(x^{\epsilon})dx^{\alpha}dx^{\beta}+
\Xi(x^{\epsilon})G_{AB}(x^C)dx^Adx^B$, $\Xi(x^{\epsilon})$ warping
factor,
\item generalized decomposable or double warped:\\ 
$ds^2=\Xi_1(x^{C})G_{\alpha\beta}(x^{\epsilon})dx^{\alpha}dx^{\beta}+\Xi_2(x^{\epsilon})G_{AB}(x^C)dx^Adx^B$,\ \-\ \-
$\Xi_1(x^C)$, $\Xi_2(x^{\epsilon})$ warping factors,
\item twisted product:\\
$ds^2=G_{\alpha\beta}(x^{\epsilon})dx^{\alpha}dx^{\beta}+\Xi_2(x^{a})G_{AB}(x^C)dx^Adx^B$
\item conformally reducible:\\
$ds^2=\Xi(x^a)(G_{\alpha\beta}(x^{\epsilon})dx^{\alpha}dx^{\beta}+G_{AB}(x^C)dx^Adx^B)$,
\item bi-conformally flat:\\
$ds^2=\Xi_1(x^a)\eta_{\alpha\beta}(x^\epsilon)dx^{\alpha}dx^{\beta}+\Xi_2(x^a)\eta_{AB}(x^C)dx^Adx^B$,
 $\eta_{\alpha\beta}$, $\eta_{AB}$ flat metrics of dimension $p$ and $n-p$ respectively.
\end{enumerate}
\label{classification}
\end{defn} 

The coordinate system of definitions \ref{separable} and \ref{classification} 
is fully adapted to the decomposition of the metric
tensor but in general we cannot expect this to be the case. 
Therefore it would be desirable to have a result characterizing 
conformally separable pseudo-Riemannian manifolds 
or any of the cases presented in definition \ref{classification}
in a coordinate-free way. 
Next we prove an intrinsic {\em local} characterization 
valid for a general conformally separable pseudo-Riemannian manifold 
which will enable us to derive characterizations for most of 
the particular 
cases described in definition \ref{classification} in a simple way. 
A lemma is needed first.    
\begin{lem}
The following assertion is true
\be
T_{abc}=0\Longleftrightarrow
\bnb_aP_{bc}=-\fr{1}{p}E_aP_{bc},\ \bnb_a\Pi_{bc}=-\fr{1}{n-p}W_a\Pi_{bc}.
\label{assertion}
\ee
\label{ftp}
\end{lem} 

\begin{pf}
First of all, it is convenient to rewrite the condition $T_{abc}=0$ in an appropriate form.
 From (\ref{tensor-t}) we have 
\be
M_{abc}=\fr{1}{p}E_aP_{bc}-\fr{1}{n-p}W_a\Pi_{bc}.
\label{star}
\ee
Using the definition $M_{abc}=\nb_bP_{ac}+\nb_cP_{ab}-\nb_aP_{bc}$ we can isolate $\nb_bP_{ac}$ getting
\be
\nb_bP_{ac}=\fr{1}{2p}(E_aP_{bc}+E_cP_{ab})-\fr{1}{2(n-p)}(W_a\Pi_{bc}+W_c\Pi_{ba}).
\label{25}
\ee
Each (\ref{star}) and (\ref{25}) are equivalent to $T_{abc}=0$. Next we show the equivalence of $T_{abc}=0$ to 
the conditions (\ref{assertion}). Expanding 
$\bnb_aP_{bc}$ by means of (\ref{identity-1}) and use of (\ref{assertion})  
yields
\bea
\nb_bP_{ac}&=&\fr{1}{2p}(E_aP_{bc}+E_cP_{ab})+\fr{1}{2}(P_{cp}M^p_{\ ba}+P_{ap}M^p_{\ bc})\\
\nb_b\Pi_{ac}&=&\fr{1}{2(n-p)}(W_a\Pi_{bc}+W_c\Pi_{ab})-\fr{1}{2}(\Pi_{cp}M^p_{\ ba}+\Pi_{ap}M^p_{\ bc}),
\eea
which are equivalent to (write $\nb_aP_{bc}$, $\nb_a\Pi_{bc}$ in terms of $M_{abc}$)
\be
P_{cp}M^p_{\ ab}=-\fr{1}{n-p}W_c\Pi_{ab},\ \ \Pi_{cp}M^p_{\ ab}=\fr{1}{p}E_cP_{ab},
\ee
whose addition leads to $T_{abc}=0$. Conversely, suppose that $T_{abc}=0$. Then inserting (\ref{star}) and 
(\ref{25}) into (\ref{identity-1}) gives us the condition on $\bnb_aP_{bc}$ at once.
 The calculation for $\bnb_a\Pi_{bc}$ is similar using the identities written in terms of $\Pi_{ab}$.\qed
\end{pf}

\begin{thm}
A pseudo-Riemannian manifold $(V,\G)$ is conformally separable at the point $p\in V$ if and only if 
there exists an orthogonal projector 
$P_{ab}$ such that the tensor $T_{abc}$ formed with $P_{ab}$ and its complementary $\Pi_{ab}=\rmg_{ab}-P_{ab}$ is 
zero identically in a neighbourhood of $p$. 
In such case $P_{ab}$ and $\Pi_{ab}$ are the leaf metrics of the separation.
\label{conf-separable}
\end{thm}

\begin{pf} 
To show that the condition of the theorem is necessary we simply choose 
the local coordinates around $p$ in which the metric tensor takes the 
form of (\ref{metric-separable}) and calculate 
the tensor $T_{abc}$ as in example \ref{example1}. 
Use of (\ref{27}) readily implies that 
$T_{abc}=0$. To prove that the condition is also sufficient
 choose an orthonormal co-basis 
$\{\bz^1,\dots,\bz^n\}$ adapted to $P_{ab}$ and $\Pi_{ab}$, that is to say, 
(we use index-free 
notation and index label splitting as in definition \ref{separable})
$$
\bP=\sum_{\alpha=1}^p\epsilon_{\alpha}\bz^{\alpha}\otimes\bz^{\alpha},\ 
\bPi=\sum_{A=p+1}^{n}\epsilon_{A}\bz^A\otimes\bz^A,
$$
where $\epsilon_{\alpha},\ \epsilon_A=\pm 1$ 
(the exact value for each index $\alpha$, $A$ 
will depend on the signature of $\rmg_{ab}$).
Now since 
$$
\bnb_c\bz^{\alpha}=
-\bbg^{\alpha}_{\ bc}\bz^b=-\bbg^{\alpha}_{\ \beta c}\bz^{\beta}-\bbg^{\alpha}_{\ Bc}\bz^B, 
$$
we have
$$
\bnb_c\bP=
-\sum_{\alpha=1}^{p}\epsilon_{\alpha}[\bbg^{\alpha}_{\ \beta c}(\bz^{\beta}\otimes\bz^{\alpha}+
\bz^{\alpha}\otimes\bz^{\beta})+
\bbg^{\alpha}_{\ Bc}(\bz^{B}\otimes\bz^{\alpha}+\bz^{\alpha}\otimes\bz^B)],
$$
which by (\ref{assertion}) entails $\bbg^{\alpha}_{\ Bc}=0$. Similarly condition 
(\ref{assertion}) on ${\bf \Pi}$ 
shows that 
$\bbg^B_{\ \alpha c}=0$. These two conditions upon the connection coefficients 
imply by means of Frobenius theorem that
the distributions spanned by $\{\bz^1,\dots,\bz^p\}$ and $\{\bz^{p+1},\dots,\bz^n\}$ 
are both integrable. 
Therefore in a local coordinate system $\{x^1,\dots,x^n\}$ around $p$ adapted to 
the manifolds generated by these distributions 
(i.e in these coordinates the manifolds are given by the conditions 
$x^{\alpha}=c^{\alpha}$, $x^A=c^A$) the metric tensor takes the form 
$$
ds^2=\rmg_{\alpha\beta}(x^a)dx^{\alpha}dx^{\beta}+\rmg_{AB}(x^a)dx^Adx^B,
$$
and the tensors $P_{ab}$ and $\Pi_{ab}$ look like
$$
P_{ab}=\rmg_{\alpha\beta}\delta^{\alpha}_{\ a}\delta^{\beta}_{\ b},\ \ \ 
\Pi_{ab}=\rmg_{AB}\delta^A_{\ a}\delta^B_{\ b},
$$  
so the non-zero components of the Christoffel symbols are
\bnr
\Chr^{\alpha}_{\ \beta\g}&=&
\fr{1}{2}\rmg^{\alpha\rho}(\d_{\beta}\rmg_{\g\rho}+\d_{\g}\rmg_{\rho\beta}-\d_{\rho}\rmg_{\beta\g}),\ 
\Chr^{\alpha}_{\ \beta A}=\fr{1}{2}\rmg^{\alpha\rho}\d_A\rmg_{\beta\rho},\\
\ \Chr^\alpha_{\ BA}&=&-\fr{1}{2}\rmg^{\alpha\rho}\d_{\rho}\rmg_{BA},
\Chr^A_{\ B\alpha}=\fr{1}{2}\rmg^{AD}\d_{\alpha}\rmg_{BD},\\ 
\Chr^A_{\ \alpha\beta}&=&-\fr{1}{2}\rmg^{AD}\d_D\rmg_{\beta\alpha},\ 
\Chr^A_{\ BC}=\fr{1}{2}\rmg^{AD}(\d_B\rmg_{CD}+\d_C\rmg_{DB}-\d_D\rmg_{BC}),
\enr
where
$$
\rmg^{\alpha\rho}\rmg_{\rho\beta}=\delta^{\alpha}_{\ \rho},\ \ 
\rmg^{AD}\rmg_{DB}=\delta^A_{\ B}.
$$
The only nonvanishing components of $M_{abc}$, $E_a$, $W_a$ are thus
\bea
M_{\alpha AB}=\d_{\alpha}\rmg_{AB},\ \ 
M_{A\alpha\beta}=-\d_{A}\rmg_{\alpha\beta},\nonumber\\
E_{A}=-\d_{A}\log|\mbox{det}(\rmg_{\alpha\beta})|,
W_{\alpha}=-\d_{\alpha}\log|\mbox{det}(\rmg_{AB})|,
\label{M-W}
\eea
from which we deduce that those of $T_{abc}$ are
\be
T_{\alpha AB}=\d_{\alpha}\rmg_{AB}+\fr{1}{n-p}\rmg_{AB}W_{\alpha},\ \ 
T_{A\alpha\beta}=\d_A\rmg_{\alpha\beta}+\fr{1}{p}\rmg_{\alpha\beta}E_A,
\label{couple}
\ee
Thus we are left with the couple of equations
(\ref{couple}) equalled to zero. The general solution of the resulting PDE system is 
$$
\rmg_{\alpha\beta}=G_{\alpha\beta}(x^{\delta})e^{\Lambda_1(x^a)},\ 
\rmg_{AB}=G_{AB}(x^D)e^{\Lambda_2(x^a)},
$$
where $G_{\alpha\beta}$, $G_{AB}$, $\Lambda_1$, $\Lambda_2$ 
are arbitrary functions of their respective arguments 
 with no restrictions other than det$(G_{\alpha\beta})\neq 0$, det$(G_{AB})\neq 0$. 
Comparing these 
expression with (\ref{metric-separable}) the result follows.\qed
\end{pf}

\begin{rem}
{\em A global 
characterization of conformally separable pseudo-Rie\-man\-ni\-an 
manifolds was first given in \cite{YANO2} and is this: 
a pseudo-Riemannian manifold is conformally separable iff there exist 
two orthogonal families of foliations by totally umbilical hypersurfaces.
The family of first fundamental forms of each hypersurface gives rise to 
the leaf metrics of the decomposition of $\rmg_{ab}$ in the obvious way
(this result was re-derived in \cite{PONGE}).} 
\end{rem}

Theorem \ref{conf-separable}
clearly states the geometric relevance of $T_{abc}$ as a tool 
to characterize conformally separable 
pseudo-Riemannian manifolds. In fact 
the condition $T_{abc}=0$ can be re-written in terms of the factors $\Xi_1$ 
and $\Xi_2$ introduced in the definition
of a conformally separable metric. To that end we use the equivalent 
condition (\ref{25}) and replace the 1-forms 
$E_a$, $W_a$ by their expressions given in (\ref{M-W}) which can be written
as
\be
E_a=-p\Pi_a^{\ r}\d_r\log|\Xi_1|,\ \ W_a=(p-n)P_a^{\ r}\d_r\log|\Xi_2|,
\label{E-W}
\ee
whence
\be
\nb_bP_{ac}=P_{bc}u_a+P_{ab}u_c-P_a^{\ r}u_r\rmg_{bc}-P_c^{\ r}u_r\rmg_{ab},
\label{Xi12}
\ee
where 
$$
u_a=\fr{E_a}{2p}+\fr{W_a}{2(n-p)}.
$$
It is not difficult now to to characterize intrinsically 
almost all the subcases presented in definition \ref{classification}.
\begin{thm}
Under the hypotheses of theorem \ref{conf-separable} a pseudo-Riemannian 
manifold $(V,\G)$ is locally
\begin{enumerate}
\item decomposable or reducible if and only if $E_a=W_a=0$,
\item a warped product if and only if $E_a=0$ and $W_a$ is an exact 1-form,
\item a double warped product if and only if both $E_a$, $W_a$ are exact 
1-forms,
\item a twisted product if and only if $E_a=0$,
\item  conformally reducible if and only if $u_a$ is exact. 
\end{enumerate}
In all cases the conditions are understood to hold only in a neighbourhood
of a point $p$.
\label{theorem-classification}
\end{thm}  
\begin{pf}
To show that the above conditions are necessary 
we only have to apply formula (\ref{E-W}) case by case and take into
account that under the conditions of theorem \ref{conf-separable}
$$
P^a_{\ b}=\left\{\begin{array}{c}
\delta^{\alpha}_{\ \beta},\ a=\alpha,\ b=\beta,\\
0\ \mbox{otherwise.}
\end{array}\right.,\ \ 
\Pi^a_{\ b}=\left\{\begin{array}{c}
\delta^A_{\ B},\ a=A, b=B\\
0\ \mbox{otherwise.}
\end{array}\right.,
$$
where the local coordinates of (\ref{metric-separable}) have been set around
$p$.
\begin{itemize}
\item $(V,\G)$ decomposable or reducible $\Rightarrow$ $\Xi_1=\Xi_2=1$ 
$\Rightarrow$ $E_a=W_a=0$. 
\item $(V,\G)$ warped product $\Rightarrow$ $E_a=0$, $W_a=(p-n)\d_a\log|\Xi_2|$.
\item $(V,\G)$ double warped $\Rightarrow$ $E_a=-p\d_a\log|\Xi_1|$, $W_a=(p-n)\d_a\log|\Xi_2|$.
\item $(V,\G)$ twisted product $\Rightarrow$ $\Xi_1=1$ $\Rightarrow$ $E_a=0$.
\item $(V,\G)$ conformally reducible $\Rightarrow$ $\Xi_1=\Xi_2=\Xi$ $\Rightarrow$
$u_a=-\d_a\log|\Xi|^{1/2}$.
\end{itemize}
The sufficiency follows from simple algebraic manipulations involving the 
conditions of each case and the relations 
$$
E_A=-p\fr{\d_A\Xi_1}{\Xi_1},\ W_A=-(n-p)\fr{\d_{\alpha}\Xi_2}{\Xi_2},
$$ 
coming from (\ref{27}) which holds due to the property $T_{abc}=0$. 
\begin{itemize}
\item $E_a=W_a=0$ $\Rightarrow$ $\d_A\Xi_1=\d_{\alpha}\Xi_2=0$ $\Rightarrow$ 
$\Xi_1=\Xi_1(x^{\alpha})$, $\Xi_2=\Xi_2(x^A)$ $\Rightarrow$ $(V,\G)$ is decomposable or 
reducible.
\item $E_a=0$ and $W_a$ exact $\Rightarrow$ $\Xi_1=\Xi_1(x^{\alpha})$ and for some 
scalar function $\Phi$ we have
$$
-(n-p)\fr{\d_{\alpha}\Xi_2}{\Xi_2}=\d_{\alpha}\Phi,\ 0=\d_A\Phi,
$$  
which entails $\Xi_2=\Xi_2(x^{\alpha})$ $\Rightarrow$ $(V,\G)$ is a warped product.
\item $E_a$ and $W_a$ are both exact. A calculation similar to the previous point 
leads to $\Xi_1=\Xi_1(x^A)$, $\Xi_2=\Xi_2(x^{\alpha})$ $\Rightarrow$ $(V,\G)$ is double warped.
\item $E_a=0$ $\Rightarrow$ $\Xi_1=\Xi_1(x^{\alpha})$ $\Rightarrow$ $(V,\G)$ is a twisted 
product.
\item If $u_a$ is exact then for some scalar function $\Phi$ we have 
$$
-\fr{1}{2}\d_{\alpha}\log|\Xi_2|=\d_{\alpha}\Phi,\ 
-\fr{1}{2}\d_A\log|\Xi_1|=\d_A\Phi
$$  
which implies that either $(V,\G)$ is a double twisted product (and in particular conformally reducible 
with $\Xi=\Xi_1(x^A)\Xi_2(x^{\alpha})$) or $|\Xi_1|=|\Xi_2|$ $\Rightarrow$ $(V,\G)$ is conformally 
reducible.\qed
\end{itemize}
\end{pf}
Local characterizations of some of the cases presented in theorem 
\ref{theorem-classification} are already known
and have been rediscovered several times by different procedures. 
For instance the reducibility condition is clearly equivalent to 
$\nb_aP_{bc}=\nb_a\Pi_{bc}=0$ which was proven in 
\cite{PETROV} in the context of General Relativity
(\cite{YANOMASAHIRO} also proves this result in Riemannian geometry). 
This result is known in Riemannian geometry as {\em de Rham decomposition theorem} 
and it was formulated by de Rham in both local and global terms \cite{DERHAM} (the 
global version was formulated as early as in 1924 \cite{BOMPIANI}). 
A local characterization of Riemannian warped products is sometimes attributed to 
Hiepko \cite{HIEPKO} but in \cite{KRUCKOVIC} such characterization is already present. 
Alternative local characterizations to those of 
theorem \ref{theorem-classification} of double 
warped products and certain conformally reducible manifolds were found in \cite{CAROT2,CAROT}  
in the framework of General relativity
and general conformally reducible Riemannian 
manifolds were locally characterized in \cite{ANCIKOV}. There are also {\em global 
characterizations} of the cases discussed in theorem \ref{theorem-classification} 
(see \cite{KUPELI} for a summary of them).

In any case our method is more general and simpler than the procedures 
followed so far and it covers virtually all possible types of 
conformally reducible pseudo-Riemannian manifolds being all of them presented in a single 
general result (theorem \ref{theorem-classification}). This makes of the bi-conformal connection
an important tool in the study of conformally separable pseudo-Riemannian manifolds and we believe
that it could play a key role in the research of these manifolds because it is the natural connection
which keeps the decomposition of the metric tensor in the two leaf metrics defined by the
projectors $P_{ab}$ and $\Pi_{ab}$. An explicit example of this key role is the local characterization 
of conformally separable pseudo-Riemannian manifolds with a conformally flat leaf metric, not covered
by theorem \ref{theorem-classification}. These cases have never been tackled before and in the companion
paper following this publication we show how the use of the bi-conformal connection allows us to 
find a local characterization in the same terms as theorem \ref{theorem-classification}. 
The impatient reader may consult all the details in \cite{GRQC} but just to give a first glimpse we
will say that this characterization comes in through the vanishing of a four rank tensor 
involving the curvature tensor of the bi-conformal connection 
(in a sense it resembles the local characterization of conformally flat 
pseudo-Riemannian manifolds by means of the vanishing of the Weyl tensor). 
More examples showing explicitly the usefulness of our techniques are 
supplied in next section.

\section{Examples}
\label{examples}
\begin{exmp}
{\em As our first example we consider the four dimensional 
pseudo-Rie\-man\-ni\-an manifold with metric given by
$$
ds^2=(\Psi^2\sin^2\z-\alpha^2)dt^2+2\Psi^2\sin^2\z d\phi dt+B^2(dr^2+r^2d\z^2)+\Phi^2\sin^2\z d\phi^2,
$$
where the coordinate ranges are $-\infty<t<\infty,\ 0<r<\infty,\ 0<\z<\pi,\ 0<\z<2\pi$ and
the functions $\Psi$, $\alpha$, $B$ and $\Phi$ only depend on the the coordinates $r$, $\z$. 
We will try to find out the conditions under which the metric 
is conformally separable with the hypersurfaces $t=const$ as one of the leaves.  
A simple calculation shows that 
the projector $P^a_{\ b}$ projecting vectors onto the distribution generated by the 
above hypersurfaces  is (now and henceforth all the components omitted in 
an explicit tensor representation are understood to be zero)
$$
P^r_{\ r}=P^{\z}_{\ \z}=P^{\phi}_{\ \phi}=1,\ P^{\phi}_{\ t}=\fr{\Psi^2}{\Phi^2},
$$
which entails
$$
P_{tt}=\fr{\Psi^4}{\Phi^2}\sin^2\z,\ P_{rr}=B^2,\ P_{\z\z}=r^2B^2,\ P_{\phi\phi}=\Phi^2\sin^2\z,\ 
P_{t\phi}=\Psi^2\sin^2\z.
$$
From here we can calculate the components of the tensor $T_{abc}$ and set them equal to zero. 
After doing this we find the following independent conditions (letter subscripts 
mean partial derivatives)
$$
-\Psi\Phi_r+\Psi_r\Phi=0,\ \ -\Psi\Phi_{\z}+\Psi_{\z}\Phi=0,
$$
which are fulfilled if and only if
$$
|\Phi|=|\Psi|
$$
Under these conditions the metric takes the form
$$
ds^2=\Psi^2\sin^2\z(dt+d\phi)^2-\alpha^2dt^2+B^2(dr^2+r^2d\z^2).
$$
This metric is not written in the form of (\ref{metric-separable}) and so it is not evident that
it is conformally separable with the hypersurfaces $t=const$ as the leaves. This is so because the 
remaining coordinates $r$, $\z$, $\phi$ are not adapted to the separation and so a coordinate 
change would be necessary to bring the above metric into the form (\ref{metric-separable}). 
An advantage of our technique is that we do not need to find this coordinate change 
and only by prescribing one of the leaves of the separation have we been able to determine easily
that our pseudo-Riemannian manifold is conformally separable. In this particular example we can even
go further and calculate the 1-forms $E_a$ and $W_a$. In this way we obtain that $E_a=0$ whereas 
$W_a$ is closed (locally exact) so theorem \ref{theorem-classification} says that this pseudo-Riemannian 
manifold is locally a warped product.}
\end{exmp}
\begin{exmp}
\label{example2}
{\em In the foregoing results we have only concentrated on conformally separable 
pseudo-Riemannian manifolds but nothing was said about manifolds with 
{\em conformal slices} and not 
conformally separable. To illustrate this case let us consider the four dimensional pseudo-Riemannian manifold 
given in local coordinates $\{x^1,x^2,x^3,x^4\}$ by
\bea
\hspace{-.8cm}ds^2=\Phi(x)[\Xi_1(x^1,x^2,x^3)(dx^1)^2+
\Xi_2(x^1,x^2,x^3)(dx^2)^2+\Xi_3(x^1,x^2,x^3)(dx^3)^2]+\nonumber\\
+2\sum_{i=1}^3\beta_i(x)dx^idx^4+\Psi(x)(dx^4)^2,\hspace{2cm}
\label{general-foliation}
\eea
where $x=\{x^1,x^2,x^3,x^4\}$ and $\Phi(x)$, $\beta_i(x)$, $\Psi(x)$, 
$\{\Xi_i(x^1,x^2,x^3)\}_{i=1,2,3}$ 
are functions at least $C^2$ in an open domain.
The above line element is the most general four dimensional metric admitting a
local foliation by three dimensional conformal hypersurfaces 
(here these are given by the condition $x^4=const$)
because according to a classical result
any three dimensional metric tensor can be written as 
the bracket term multiplying $\Phi(x)$ in equation (\ref{general-foliation}).
  
The non-zero components of the orthogonal projector $P^{a}_{\ b}$ associated to the 
foliation $x^4=const$ (see previous example) are
\be
P^1_{\ 1}=P^2_{\ 2}=P^3_{\ 3}=1,\ P^i_{\ 4}=\fr{\beta_i(x)}{\Phi(x)\Xi_i(x^1,x^2,x^3)},\ i=1,2,3,
\label{projection}
\ee
from which we easily get
$$
P_{11}=P_{22}=P_{33}=\Phi(x),\ P_{i4}=\beta_i(x),\ i=1,2,3,\ 
P_{44}=\sum_{i=1}^3\fr{\beta^2_i(x)}{\Phi(x)\Xi_i(x^1,x^2,x^3)}.
$$
Using this we can check the condition $T_{abc}=0$ and find out what 
is obtained. 
This is a rather long calculation which is easily 
performed with any of the computer algebra systems available today (the system used here was GRTensorII).
The result is that the tensor $T_{abc}$ does not vanish in this case although
a calculation using (\ref{projection}) shows the important property
\be
P^r_{\ a}P^s_{\ b}P^q_{\ c}T_{rsq}=0.
\label{changed-condition}
\ee
\begin{thm}
A necessary condition that a four dimensional pseudo-Rie\-man\-ni\-an 
manifold can be foliated by 
conformal hypersurfaces with associated orthogonal projector $P_{ab}$ is equation 
(\ref{changed-condition}).
\end{thm}  
 This result suggests 
that it may well be possible to generalize conditions of theorem \ref{conf-separable} 
to metrics of arbitrary 
dimension which are not conformally separable 
replacing these conditions by (\ref{changed-condition}). 
The true extent of this assertion is under current research.} 
\end{exmp}

\section*{Acknowledgements}
We wish to thank Jos\'e M. M. Senovilla for a careful reading of 
the manuscript and his many suggested improvements. 
Financial support from projects\\  9/UPV00172.310-14456/2002 
and BFM 2000-0018 is also gratefully acknowledged. 
Finally we thank an anonymous referee for his valuable 
comments.

\appendix
\setcounter{section}{0}
\section{Basic identities involving the Lie derivative}
In this appendix we recall some properties of the Lie derivative needed in the main text. 
Despite their basic character, they are hardly presented in basic Differential Geometry 
textbooks and the author is only aware of \cite{YANO,SCHOUTEN} as the only references in 
which they are studied.  
\begin{prop} 
For any symmetric connection $\bnb$ defined in a differentiable manifold $V$, any vector field 
$\xiv$ at least $C^2$ and a tensor field $T^{a_1\dots a_p}_{\ b_1\dots b_q}\in T^p_{q}(V)$ 
we have the following identities
\bea
\lie\bar{\g}^a_{bc}=\bnb_b\bnb_c\xi^a+\xi^d\bar{R}^a_{\ cdb},\label{lie-xi}\\
\hspace{-1cm}\bnb_c\lie T^{a_1\dots a_s}_{\ b_1\dots b_q}-
\lie\bnb_cT^{a_1\dots a_s}_{\ b_1\dots b_q}=\nonumber\\
=-\sum_{j=1}^s(\lie\bar{\g}^{a_j}_{cr})T^{\dots a_{j-1}ra_{j+1}\dots}_{\ b_1\dots b_q}+
\sum_{j=1}^q(\lie\bar{\g}^r_{cb_j})T^{a_1\dots a_s}_{\dots b_{j-1}rb_{j+1}\dots},
\label{lie-conmmutation}\\
\lie \bar{R}^d_{\ cab}=\bnb_a(\lie\bar{\g}^d_{\ bc})-\bnb_b(\lie\bar{\g}^d_{ac})
\label{lie-curvature},
\eea
where $\bar{\g}^a_{\ bc}$ are the components of the connection $\bnb$ and $\bar{R}^a_{\ bcd}$ its
curvature (these identities are calculated under
the convention (\ref{convention}) for the curvature tensor).
 Furthermore if a metric tensor $\rmg_{ab}$ is set in $V$ and $\nb$ 
is now the metric connection associated to it then
\be
\lie\g^a_{bc}=\fr{1}{2}\rmg^{ae}\left[\nabla_b(\lie\rmg_{ce})+
\nabla_c(\lie\rmg_{be})-\nabla_e(\lie\rmg_{bc})\right].\label{lie-connection}
\ee
\end{prop} 

\begin{rem}{\em
The Lie derivative of the connection is a tensor even though 
$\g^a_{\ bc}$ is not. To see this we denote by 
$\{\Phi_s\}$ the one-parameter group of local diffeomorphisms 
generated by the vector field 
$\xiv$ and by $(\Phi^*_s\g)^a_{\ bc}$ the transformed of the connection
under $\Phi_s$ which in the local coordinates $x=\{x^1,\dots,x^n\}$ 
is calculated by means of the formula
$$
(\Phi^*_s\Gamma)^a_{\ bc}(x)=
\fr{\d\Phi_{-s}^{a}}{\d x^r}\fr{\d\Phi^z_s}{\d x^b}\fr{\d\Phi^q_s}{\d x^c}\Gamma^r_{\ zq}(\Phi_s(x))+
\fr{\d\Phi_{-s}^{a}}{\d x^r}\fr{\d}{\d x^b}\left(\fr{\d\Phi_s^r}{\d x^c}\right).
$$
Hence neither $(\Phi^*_s\g)^a_{\ bc}$ nor $\g^a_{\ bc}$ are tensors 
but the difference $(\Phi^*_s\g)^a_{\ bc}-\g^a_{\ bc}$ is a tensor and 
this implies that 
$$
\lim_{s\rightarrow 0}\fr{(\Phi^*_s\g)^a_{\ bc}-\g^a_{\ bc}}{s}
$$
is also a tensor which is the Lie derivative of the connection.
}
\end{rem}





\begin{thebibliography}{}
%
%
\bibitem{ANCIKOV} An\v cikov A. M., Izv. Vys\v s. U\v cebn. Zaved. Matematika \textbf{48}, (1965) 13-18.
(English translation: Amer. Math. Soc. Transl. \textbf{92}, (1970) 201-207.)  
\bibitem{BOMPIANI} Bompiani E., Rendiconti del Circolo Matematico di Palermo \textbf{48}, (1924) 121-134.
\bibitem{CAROT} Carot J. and Tupper B. O. J., Class. Quantum Grav. \textbf{19}, (2002) 4141-4166.
\bibitem{CAROT2} Ramos M. P. M., Vaz E. G. L. R. and Carot J., J. of Math. Phys. {\bf 44}, (2003)
 4839-4865. 
\bibitem{SERGI} Coll B., Hildebrandt S. R. and Senovilla J. M. M.,
 Gen. Rel. Grav. \textbf{33}, (2001) 649-670.
\bibitem{DERHAM} de Rham G., Comment. Math. Helv. {\bf 26}, (1952) 328-344.
\bibitem{EISENHART} Eisenhart L. P., \textit{Riemannian Geometry} (Princeton University Press, Princeton 1964) p. 33.
\bibitem{EISENHARTII} Eisenhart L. P., \textit{Continuous groups of transformations} (Dover publications Inc., 
New York 1933).
\bibitem{KUPELI} Fern\'andez-L\'opez M., Garc\'{\i}a-R\'{\i}o E., Kupeli D.N.
and \"Unal B., manuscripta math. {\bf 106}, (2001) 213-217.
\bibitem{BI-CONFORMAL} Garc\'{\i}a-Parrado A. and Senovilla J.M.M., Class. Quantum Grav. \textbf{21}, (2004) 
2153-2177.
\bibitem{GRQC} Garc\'{\i}a-Parrado A., {\em Preprint} \verb|http://www.arxiv.org/abs/math-ph/0409037|
\bibitem{LIBROHALL} Hall G. S., \textit{Symmetries and Curvature Structure in General Relativity} (World Scientific, Singapore 2004).
\bibitem{HIEPKO} Hiepko S., Math. Ann. {\bf 241}, (1979) 209-215.
\bibitem{KRUCKOVIC} Kru\v ckovi\v c G., I. Trudy Sem. Vektor. Tenzor. Anal. {\bf 11}, (1961) 103-128;
erratum ibid. {\bf 12}, (1961) 95.
\bibitem{PETROV} Petrov A. Z., \textit{New methods in general relativity} (Nauka, Moscow 1966).
\bibitem{PONGE} Ponge R. and Reckziegel H., Geom. Dedicata \textbf{48}, (1993) 15-25.  
\bibitem{REINHART} Reinhart B. L., 
\textit{Differential geometry of Foliations, the fundamental integrality problem} (Springer, Berlin 1983).  
\bibitem{SCHOUTEN}
Schouten J. A., \textit{Ricci Calculus} (Springer, Berlin 1954).
\bibitem{WALKER1} Walker A. G., {\it Quart. J. Math. Oxford} {\bf 6}, (1955)
 301-308.
\bibitem{WALKER2} Walker A. G., {\it Quart. J. Math. Oxford} {\bf 9}, (1958)
221-231.
\bibitem{WILLMORE1} Willmore T. J., {\it Proc London Math. Soc.} {\bf 6}, (1956)
191-204.
\bibitem{WILLMORE2} Willmore T. J., {\it Quart. J. Math. Oxford} {\bf 7}, (1956) 
269-276.
\bibitem{YANO2} Yano K., Proc. Imp. Acad. Tokyo \textbf{16}, (1940) 83-86.
\bibitem{YANO}
Yano K., \textit{Theory of Lie Derivatives} (North Holland, Amsterdam 1955). 
\bibitem{YANOMASAHIRO} Yano K. and Masahiro K., \textit{Structures on manifolds} (World Scientific, Singapore 1984).
\end{thebibliography}
\end{document}